\documentclass[11pt]{article}
\usepackage[final]{epsfig}
\usepackage[left=3cm,right=3cm, top=2.5cm,bottom=2.5cm,bindingoffset=0cm]{geometry}
\usepackage{graphics}
\usepackage{amsmath}
\usepackage{amsfonts}
\usepackage{latexsym}
\usepackage{amssymb, mathabx}
\usepackage{amsthm}
\usepackage{graphicx}
\usepackage{epstopdf}
\usepackage{multicol,multirow}
\usepackage{hyperref, enumitem}
\usepackage{mathtools}

\usepackage{mathrsfs}

 \usepackage{relsize}

\usepackage[nottoc]{tocbibind}

\usepackage{amsthm}

\usepackage[latin1]{inputenc}
\usepackage{tikz}
\usetikzlibrary{shapes,arrows,decorations.pathmorphing}
\usepackage{float}

\usepackage{blkarray}
\usepackage{xypic}

\usepackage[title]{appendix}

\let\underbrace\LaTeXunderbrace

\numberwithin{equation}{section}

\makeatletter
\def\thmhead@plain#1#2#3{%
  \thmname{#1}\thmnumber{\@ifnotempty{#1}{ }\@upn{#2}}%
  \thmnote{ {\the\thm@notefont#3}}}
\let\thmhead\thmhead@plain

\makeatother

\newcounter{AppCounter}

\def\restrict#1{\raise-.5ex\hbox{\ensuremath|}_{#1}}

\newtheorem{lemma}{Lemma}[section]
\newtheorem{proposition}[lemma]{Proposition}

\newtheorem{remark-definition}[lemma]{Remark-Definition}
\newtheorem{theorem}[lemma]{Theorem}

\newtheorem{proposition-conjecture}[lemma]{Proposition-conjecture}
\newtheorem{problem}[lemma]{Problem}
\newtheorem{question}[lemma]{Question}

\theoremstyle{definition}

\newtheorem{remark}[lemma]{Remark}

\newcommand{\marginnote}[1]
{
}

\newcounter{cy}

\newcounter{bk}

\newcounter{dps}

\usepackage{lipsum}

\title{On the Finite-Time Singularities of Point-Vortex Systems I: The Half-Plane}

\author{ Cheng Yang\thanks{Division of Mathematical Sciences, Nanyang Technological University, 21 Nanyang Link, Singapore 637371; 
 e-mail: \tt{cheng.yang@ntu.edu.sg} }
\\
}

\begin{document}
\date{}
\maketitle

\begin{abstract}

In this paper, the first in a series devoted to finite-time singularities of point-vortex systems in domains with boundaries, we investigate finite-time boundary collisions and interior collapse for three vortices in the upper half-plane $\mathbb{H}^+$. We first prove that finite-time boundary collisions are impossible for arbitrary nonzero vortex circulations, without any restrictions on their signs or algebraic relations. The result excludes both simultaneous boundary collapse and individual vortex approaches to the boundary, without requiring convergence of the vortex trajectories. In contrast, we establish the existence of finite-time interior collapse by treating the boundary-induced interactions as a perturbation of the classical self-similar Gr\"{o}bli collapse in the full plane. Under suitable conditions on the vortex circulations and collapse profile, we prove the existence of a one-dimensional local strong stable manifold of the augmented rescaled system, which generates exact solutions in which three distinct vortices collapse simultaneously at any prescribed interior point $z^*\in\mathbb{H}^+$. Moreover, these solutions are asymptotically self-similar, with the characteristic square-root collapse rate $O(\sqrt{T-t})$.

\end{abstract}

\tableofcontents

\section{Introduction} \label{intro}

\subsection{Classical Point-Vortex Dynamics and Finite-Time Collapse}

The point-vortex model, originating from Helmholtz's seminal
work \cite{Helmh1858}, provides a finite-dimensional description
of ideal incompressible planar flows with concentrated vorticity.
In this model, the vorticity distribution is represented by
\[
\omega(x,t)=\sum_{k=1}^N\Gamma_k\delta(x-x_k(t)),
\qquad
\Gamma_k\in\mathbb R\setminus\{0\},
\]
where $x_k(t)\in\mathbb R^2$ denotes the position of the
$k$th vortex and $\Gamma_k$ its circulation.
The vortex centers evolve according to a nonlinear Hamiltonian
system, whose classical formulation was established by
Kirchhoff \cite{Kirch1876}. In the full plane $\mathbb R^2$,
the conservation laws associated with translational and
rotational symmetries imply integrability for $N\leq 3$
\cite{Grobli1877,Poinc1893,Aref1979}, whereas systems with
$N\geq4$ are generally non-integrable \cite{Ziglin1980,Aref1983}.

An important feature of point-vortex dynamics is the
possibility of finite-time singularity formation, in which
two or more initially distinct vortices approach a common
point. The classical example is the self-similar collapse
of three vortices discovered by Gr\"{o}bli
\cite{Grobli1877} and subsequently studied by
Synge \cite{Synge1949}, Novikov \cite{Novikov1975},
and Aref \cite{Aref1979}.
For such a collapse, the circulations necessarily satisfy
$\Gamma_1\Gamma_2+\Gamma_1\Gamma_3+\Gamma_2\Gamma_3=0,$
and, for suitable initial configurations, the vortex
triangle preserves its shape while contracting and rotating,
with a characteristic length scale proportional to
$\sqrt{T-t}$ as the collapse time $T$ is approached.
The existence and structure of self-similar collapse for
larger vortex systems were further investigated by
O'Neil \cite{ONeil1987}. Although such singularities are
possible, they are nongeneric: under suitable assumptions,
the set of initial configurations leading to finite-time
collision has Lebesgue measure zero \cite{MP1994}.

\subsection{Finite-Time Singularities in the Half-Plane}

The presence of a solid boundary introduces additional
interactions and substantially changes the structure of
point-vortex dynamics. The Hamiltonian formulation in
general domains was developed by Lin \cite{Lin1941,Lin1943}
using the Green's function and its associated Robin function.

In the upper half-plane
$\mathbb H^+=\{z\in\mathbb C:\operatorname{Im}z>0\},$
the loss of vertical translational and rotational symmetry leaves horizontal translation as the only continuous spatial symmetry. 
Consequently, the classical integrability
of the three-vortex problem in $\mathbb R^2$ does not
extend directly to the half-plane.

These boundary-induced interactions also distinguish
two possible types of finite-time singularities:
interior collapse, in which several vortices converge to
a common point in $\mathbb H^+$, and boundary collision,
in which one or more vortices approach
$\partial\mathbb H^+$.
In our previous work \cite{Yang}, we proved that two
point vortices in the half-plane can neither collide
with each other nor reach the boundary in finite time.
It is therefore natural to investigate whether finite-time
singularities can occur for three or more vortices.
We ask the following question in \cite{Yang}.
\begin{problem}
Is it possible that an $N$-vortex system in the half-plane collides or hits the
boundary in finite time for $N>2$?
\end{problem}
In the present paper, we address this fundamental open problem concerning finite-time singularities of point vortices in the upper half-plane $\mathbb{H}^+$. The corresponding problem in the quadrant will be studied in a forthcoming paper \cite{Yang2}.

Newton \cite{Newton2001} and Aref \cite{Aref2007} note that many questions remain open about whether point vortices in a domain with a boundary can collapse or collide with the boundary. Newton explicitly highlights the problem in his book \cite{Newton2001} as follows.
\begin{quote}
``\dots a general study of the $N$-vortex collision problem has not yet been carried out. For instance, we do not know if there are any non self-similar collisions of vortices in the plane or on the sphere, and vortex collisions in domains with boundaries are only beginning to be investigated. In general, much more is known about $N$-body collisions in the context of celestial mechanics, and some of these techniques might prove useful.''
\end{quote}

\subsection{Main Results and Organization of the Paper}

This paper is the first in a series on singularity formation in point-vortex systems in domains with boundaries. We prove that three point vortices in the upper half-plane $\mathbb{H}^+$ can collapse in finite time at an interior point, and we rule out finite-time collisions with the boundary. Specifically, we establish the following
main results.

\paragraph{Finite-time boundary non-collision for three vortices}
We first establish that finite-time boundary collisions are impossible for three point vortices in the upper half-plane, for arbitrary nonzero circulations, without any restrictions on their signs or algebraic relations.  More precisely, for every finite $T>0$ such that the
solution exists on $[0,T)$
\[
\inf_{0\le t<T}\operatorname{Im} z_i(t)>0,
\qquad i=1,2,3.
\]
In particular, our result excludes both simultaneous boundary collapse and the approach of any individual vortex to the boundary, without assuming that the vortex positions or their distances from the boundary converge as \(t\uparrow T\).

The proof relies on conservation of the Hamiltonian and vertical impulse, together with suitable geometric estimates and differential inequalities. When the vertical impulse is nonzero, conservation laws and Gr\"{o}nwall-type estimates exclude boundary collisions. When the vertical impulse vanishes, energy conservation forces the two vortices with circulations of the same sign to form a tightly concentrated pair near the boundary, then appropriate monotonicity arguments rule out finite-time boundary approach. The analysis distinguishes between zero and nonzero total circulation, for which different monotonicity mechanisms are required.

\paragraph{Existence of finite-time interior collapse for
three vortices.}
We prove that finite-time interior collapse is possible in
the upper half-plane, despite the presence of boundary-induced
interactions. For three nonzero vortex strengths satisfying
$\Gamma_1\Gamma_2+\Gamma_1\Gamma_3+\Gamma_2\Gamma_3=0,$
and an admissible non-equilateral Gr\"{o}bli collapse profile,
we establish the existence of a one-dimensional local strong
stable manifold of the augmented rescaled vortex system.
Our approach takes the classical Gr\"{o}bli collapse in the full
plane as a reference solution and treats the boundary-induced
interactions as a perturbation. In suitable similarity
variables, this perturbation vanishes as the collapse time
is approached, allowing us to construct exact collapsing
solutions through a strong stable manifold argument.

In particular, we show that three distinct vortices can
collapse simultaneously at any prescribed interior point
$C_*\in\mathbb{H}^+$ in finite time $T$. The resulting
trajectories satisfy the asymptotic expansion
\[
z_k(t)
=C_*+\sqrt{T-t}\,e^{-i\Omega s(t)}W_k^*
+O(T-t),
\qquad
s(t)=-\log(T-t),
\]
where $\mathbf W^*=(W_1^*,W_2^*,W_3^*)$ denotes the
corresponding Gr\"{o}bli collapse profile. Thus, the collapse
is asymptotically self-similar, with the characteristic
square-root scaling. 

Together, these results characterize the possible locations of finite-time collapse for three vortices in the upper half-plane: interior collapse can occur, whereas boundary collisions are excluded.

The remainder of the paper is organized as follows. Section~\ref{sec:pre} presents the governing equations for $N$ point vortices in the upper half-plane, together with their Hamiltonian and complex formulations. We also recall the classical finite-time self-similar Gr\"{o}bli collapse of three vortices in the full plane, which serves as the reference solution for our perturbative construction. In Section~\ref{sec:bd_noncolli}, we establish the impossibility of finite-time boundary collisions for three vortices with arbitrary nonzero circulations (Theorem~\ref{thm:main_bd_noncolli}). Finally, Section~\ref{sec:rigorous_existence_collapse} develops the rescaled dynamical framework and proves the existence of finite-time interior collapsing solutions through a strong stable manifold argument (Theorem~\ref{thm:interior_collapse}).




\medskip


\section{Preliminaries}\label{sec:pre}

\subsection{Point-Vortex Dynamics in the Upper Half-Plane}
\label{sec:half_plane_dynamics}

Consider $N$ point vortices with nonzero circulations
$\Gamma_i$ and positions $z_i=x_i+iy_i$ in the upper
half-plane
\[
\mathbb H^+=\{z\in\mathbb C:\operatorname{Im}z>0\}.
\]
By the method of images, the corresponding
Kirchhoff--Routh Hamiltonian is
\begin{equation}\label{eq:Hamiltonian_func}
\begin{aligned}
H_{\mathbb H^+}
={}&\frac{1}{4\pi}\sum_{1\le i<j\le N}
\Gamma_i\Gamma_j
\log\frac{(x_i-x_j)^2+(y_i+y_j)^2}
{(x_i-x_j)^2+(y_i-y_j)^2}\\
&+\frac{1}{4\pi}\sum_{i=1}^N
\Gamma_i^2\log(2y_i).
\end{aligned}
\end{equation}
The equations of motion take the Hamiltonian form
\begin{equation}\label{eq:Hamiltonian_eq}
\Gamma_i\dot x_i=\frac{\partial H_{\mathbb H^+}}
{\partial y_i},
\qquad
\Gamma_i\dot y_i=-\frac{\partial H_{\mathbb H^+}}
{\partial x_i},
\qquad i=1,\ldots,N.
\end{equation}

For convenience, set
\[
R_{ij}^{\pm}
=(x_i-x_j)^2+(y_i\pm y_j)^2.
\]
Differentiating \eqref{eq:Hamiltonian_func} yields
\begin{equation}\label{eq:half_plane}
\begin{aligned}
\dot x_i
={}&\frac{\Gamma_i}{4\pi y_i}
+\frac{1}{2\pi}\sum_{j\ne i}\Gamma_j
\left(
\frac{y_i+y_j}{R_{ij}^+}
-\frac{y_i-y_j}{R_{ij}^-}
\right),\\
\dot y_i
={}&\frac{1}{2\pi}\sum_{j\ne i}\Gamma_j
(x_i-x_j)
\left(
\frac{1}{R_{ij}^-}
-\frac{1}{R_{ij}^+}
\right).
\end{aligned}
\end{equation}
Here, the term $\Gamma_i/(4\pi y_i)$ represents
the horizontal drift induced by the vortex's own
image. The terms involving $R_{ij}^-$ describe
direct vortex interactions, whereas those involving
$R_{ij}^+$ arise from image vortices with
circulations $-\Gamma_j$ reflected across the
boundary $y=0$.

Equivalently, in complex notation,
\eqref{eq:half_plane} becomes
\begin{equation}\label{eq:complex}
\dot{\bar z}_k
=\frac{1}{2\pi i}\sum_{j\ne k}
\frac{\Gamma_j}{z_k-z_j}
-\frac{1}{2\pi i}\sum_{j=1}^N
\frac{\Gamma_j}{z_k-\bar z_j},
\qquad k=1,\ldots,N.
\end{equation}
The $j=k$ term in the second sum accounts for
the self-induced horizontal drift.

For two vortices in the upper half-plane, it
was proved in \cite{Yang} that neither mutual
collision nor collision with the boundary can
occur in finite time.

\subsection{Classical Gr\"obli Collapse in the Plane}
\label{sec:gna_collapse}

In contrast to the two-vortex noncollision result,
three vortices in the unbounded plane can undergo
finite-time self-similar collapse. This classical
phenomenon was first studied by Gr\"obli
\cite{Grobli1877} and subsequently revisited by
Novikov \cite{Novikov1975} and Aref \cite{Aref1979}.

For three vortices in $\mathbb C$, the equations
of motion are
\begin{equation}\label{eq:3_complex}
\dot{\bar z}_k
=\frac{1}{2\pi i}\sum_{j\ne k}
\frac{\Gamma_j}{z_k-z_j},
\qquad k=1,2,3.
\end{equation}
Conservation of the planar Hamiltonian
\[
H=-\frac{1}{2\pi}\sum_{i<j}
\Gamma_i\Gamma_j\log|z_i-z_j|
\]
implies that self-similar collapse requires
\begin{equation}\label{eq:zero_energy_cond}
\sum_{i<j}\Gamma_i\Gamma_j=0
\quad\Longleftrightarrow\quad
\frac1{\Gamma_1}+\frac1{\Gamma_2}
+\frac1{\Gamma_3}=0.
\end{equation}
Conservation of angular impulse further requires
\begin{equation}\label{eq:zero_angular_impulse}
\sum_{k=1}^3\Gamma_k|z_k(0)-z_*|^2=0,
\end{equation}
where
\[
z_*=\frac{\sum_{k=1}^3\Gamma_k z_k(0)}
{\sum_{k=1}^3\Gamma_k}
\]
is the center of vorticity and the collapse point.

More precisely, writing $Z_k=z_k(0)-z_*$, a
self-similar collapse exists if and only if
there is a constant $C\in\mathbb C$ with
$\operatorname{Re}C<0$ such that
\begin{equation}\label{eq:self_similar_shape}
\frac{1}{2\pi i}\sum_{j\ne k}
\frac{\Gamma_j}{Z_k-Z_j}
=C\overline{Z_k},
\qquad k=1,2,3.
\end{equation}
Indeed, the ansatz $z_k(t)=z_*+Z_k f(t)$
reduces \eqref{eq:3_complex} to
$f\overline{\dot f}=C$, giving
\begin{equation}\label{eq:collapse_trajectory}
z_k(t)=z_*+Z_k\sqrt{1-\frac{t}{T}}\,
\exp\left[-i\omega\log\left(1-\frac{t}{T}\right)\right],
\end{equation}
where
\[
T=-\frac{1}{2\operatorname{Re}C}>0,
\qquad
\omega=-T\operatorname{Im}C.
\]
Thus, all inter-vortex distances vanish at the
rate $(T-t)^{1/2}$. If $\omega\ne0$, the
configuration undergoes infinitely many rotations
as $t\nearrow T$.

\begin{remark}[(Role of circulation signs)]
For the point-vortex system in the upper half-plane,
solutions with circulations of the same sign exist
globally in time, with neither mutual collisions nor
boundary collisions occurring. This follows from
conservation of the Hamiltonian and the vertical
impulse $I=\sum_{i=1}^N\Gamma_i y_i$.

For mixed-sign circulations, these conservation laws
alone no longer rule out finite-time singularities.
Potential singularity scenarios include collisions
with the boundary and mutual collisions in the
interior. Boundary collisions have been investigated
under additional assumptions; see, for instance,
\cite{DoGoIf}.

In this paper, we establish a finite-time boundary
noncollision result and demonstrate the existence
of finite-time interior collisions for three vortices
in the upper half-plane.
\end{remark}


\medskip

\section{Finite-time boundary non-collision for three point vortices}\label{sec:bd_noncolli}

Let
\[
\mathbb H^+=\{z\in\mathbb C:\operatorname{Im}z>0\}
\]
denote the upper half-plane. We consider three point
vortices with positions
\[
z_i(t)=x_i(t)+iy_i(t)\in\mathbb H^+,
\qquad i=1,2,3,
\]
and nonzero real circulations
$\Gamma_1,\Gamma_2,\Gamma_3$.

Our main result excludes finite-time boundary
collisions without imposing any sign or algebraic
condition on the circulations.

\begin{theorem}[(Boundary non-collision)]
\label{thm:main_bd_noncolli}
Let $(z_1,z_2,z_3)$ be a classical collision-free
solution of the standard point-vortex equations
in $\mathbb H^+$, with arbitrary nonzero real
circulations $\Gamma_1,\Gamma_2,\Gamma_3$.

Then, for every finite $T>0$ such that the
solution exists on $[0,T)$,
\[
\inf_{0\le t<T}y_i(t)>0,
\qquad i=1,2,3.
\]
In particular, no vortex can approach the
boundary in finite time, even in the sense
of a vanishing lower limit.
\end{theorem}
The proof is divided into three cases according
to the conserved vertical impulse and the total
circulation. We first establish several estimates
that will be used throughout.

\subsection{Equations and conserved quantities}

For $i\neq j$, set
\[
R_{ij}^2=|z_i-z_j|^2,
\qquad
S_{ij}^2=|z_i-\overline z_j|^2.
\]
The half-plane point-vortex equations are
\begin{equation}
\label{eq:system}
\begin{aligned}
\dot x_i
={}&\frac{\Gamma_i}{4\pi y_i}
+\frac1{2\pi}\sum_{j\ne i}\Gamma_j
\left(
\frac{y_i+y_j}{S_{ij}^2}
-\frac{y_i-y_j}{R_{ij}^2}
\right),\\
\dot y_i
={}&\frac2\pi\sum_{j\ne i}\Gamma_j
\frac{(x_i-x_j)y_i y_j}
{R_{ij}^2S_{ij}^2}.
\end{aligned}
\end{equation}
The Hamiltonian is $H=E/(4\pi)$, where
\begin{equation}
\label{eq:energy}
E=
\sum_{i=1}^3\Gamma_i^2\log(2y_i)
+\sum_{i<j}\Gamma_i\Gamma_j
\log\frac{S_{ij}^2}{R_{ij}^2}.
\end{equation}
Consequently,
\begin{equation}
\label{eq:Econserved}
E(t)=E(0)=:E_0.
\end{equation}
The vertical impulse
\begin{equation}
\label{eq:impulse}
P(t)=\sum_{i=1}^3\Gamma_i y_i(t)
\end{equation}
is also conserved. This follows directly from
the second equation in \eqref{eq:system}, because
the contributions of each pair $(i,j)$ and $(j,i)$
cancel after multiplication by $\Gamma_i$.
Thus
\begin{equation}
\label{eq:Pconserved}
P(t)=P(0).
\end{equation}
We will repeatedly use the following elementary
barrier principle.

\begin{lemma}[(Conditional differential barriers)]
\label{lem:barrier}
Let $f:[0,T)\to(0,\infty)$ be continuously
differentiable, where $T<\infty$, and let
$\varepsilon>0$.

\textup{(i)} If
\[
f'(t)\ge-Cf(t)
\quad\text{whenever }f(t)\le\varepsilon,
\]
then
\[
f(t)\ge
\min\{f(0),\varepsilon\}e^{-Ct}.
\]

\textup{(ii)} If
\[
f'(t)\ge-Cf(t)^2
\quad\text{whenever }f(t)\le\varepsilon,
\]
then
\[
\frac1{f(t)}
\le
\max\left\{\frac1{f(0)},\frac1\varepsilon\right\}
+Ct.
\]
Both conclusions remain valid when the
differential inequalities hold only on a
subinterval, with its left endpoint replacing
the initial time.
\end{lemma}

\begin{proof}
For part (i), define
\[
F(t)=
\left(\log\frac{\varepsilon}{f(t)}\right)_+.
\]
The function $F$ is locally absolutely continuous.
Whenever $f<\varepsilon$,
\[
F'=-\frac{f'}f\le C,
\]
and whenever $f>\varepsilon$, one has $F'=0$.
Hence $F'\le C$ almost everywhere, and
\[
F(t)\le F(0)+Ct.
\]
The claimed lower bound follows.

For part (ii), define
\[
G(t)=\max\left\{\frac1{f(t)},\frac1\varepsilon\right\}.
\]
Whenever $f<\varepsilon$,
\[
G'=-\frac{f'}{f^2}\le C,
\]
while $G'=0$ whenever $f>\varepsilon$.
Integration yields the required estimate.
\end{proof}

\begin{lemma}[(Isolated boundary approach)]
\label{lem:isolated}
Suppose that for some $i\in\{1,2,3\}$ there
exist $\varepsilon,m>0$, with
$\varepsilon\le m/2$, such that
\[
y_i(t)\le\varepsilon
\quad\Longrightarrow\quad
y_j(t)\ge m
\quad\text{for every }j\ne i.
\]
Then
\[
\inf_{0\le t<T}y_i(t)>0
\]
for every finite interval of classical existence
$[0,T)$.
\end{lemma}

\begin{proof}
Whenever $y_i\le\varepsilon$, for $j\ne i$,
\[
R_{ij}\ge y_j-y_i\ge\frac m2,
\qquad
S_{ij}\ge y_j\ge m.
\]
Moreover,
\[
|x_i-x_j|\le R_{ij},
\qquad
y_j\le S_{ij}.
\]
Therefore
\[
\frac{|x_i-x_j|y_j}{R_{ij}^2S_{ij}^2}
\le\frac1{R_{ij}S_{ij}}
\le\frac2{m^2}.
\]
Using \eqref{eq:system},
\[
|\dot y_i|
\le
\frac4{\pi m^2}
\left(\sum_{j\ne i}|\Gamma_j|\right)y_i
=:Cy_i
\]
whenever $y_i\le\varepsilon$.

Lemma~\ref{lem:barrier}(i) now gives
\[
y_i(t)\ge
\min\{y_i(0),\varepsilon\}e^{-Ct}.
\]
Since $T<\infty$, the conclusion follows.
\end{proof}

\subsection{Nonzero vertical impulse}

We first treat the case $P(0)\ne0$.
No assumption on the total circulation is needed.
\begin{proposition}
\label{prop:Pnonzero}
Suppose
\[
P(0)=\sum_{i=1}^3\Gamma_i y_i(0)\ne0.
\]
Then, for every finite interval of classical
collision-free existence $[0,T)$,
\[
\inf_{0\le t<T}y_i(t)>0,
\qquad i=1,2,3.
\]
\end{proposition}

\begin{proof}
We first show that no pair of vortex heights
can become simultaneously arbitrarily small.

Fix distinct indices $i,j$ and let $k$ be
the remaining index. Suppose, by contradiction,
that there exists a sequence $t_n\uparrow T$
such that
\[
y_i(t_n)\to0,
\qquad
y_j(t_n)\to0.
\]
By conservation of vertical impulse,
\[
\Gamma_k y_k(t_n)
=P(0)-\Gamma_i y_i(t_n)-\Gamma_j y_j(t_n).
\]
Thus
\[
y_k(t_n)\longrightarrow\frac{P(0)}{\Gamma_k}.
\]
If $P(0)/\Gamma_k<0$, this is impossible
because $y_k(t)>0$. We therefore assume
\[
\frac{P(0)}{\Gamma_k}>0.
\]
We distinguish two cases.

\medskip
\noindent
\emph{Case 1: $\Gamma_i\Gamma_j<0$.}

Write
\[
L_{pq}
=\log\left(
1+\frac{4y_py_q}{R_{pq}^2}
\right)\ge0.
\]
The interaction of vortices $i,j$ satisfies
\[
\Gamma_i\Gamma_jL_{ij}\le0.
\]
Since $y_k(t_n)$ converges to a positive
finite limit, and Since \(y_i,y_j\to0\), we have
\[
R_{ik}(t_n),R_{jk}(t_n)\ge c_0>0
\]
for sufficiently large $n$.
Consequently,
\[
L_{ik}(t_n)=O(y_i(t_n)),
\qquad
L_{jk}(t_n)=O(y_j(t_n)).
\]
The energy identity therefore yields
\[
E_0\le
\Gamma_i^2\log(2y_i(t_n))
+\Gamma_j^2\log(2y_j(t_n))
+O(1).
\]
The right-hand side tends to $-\infty$,
contradicting conservation of energy.

Therefore two vortices with opposite circulations cannot simultaneously approach the boundary.

\medskip
\noindent
\emph{Case 2: $\Gamma_i\Gamma_j>0$.}

Let
\[
s=\operatorname{sgn}\Gamma_i
=\operatorname{sgn}\Gamma_j
\]
and define
\[
B(t)=|\Gamma_i|y_i(t)+|\Gamma_j|y_j(t).
\]
Then
\[
P(0)=sB(t)+\Gamma_k y_k(t).
\]
Since $P(0)/\Gamma_k>0$, there exist
$\varepsilon,m>0$ such that
\[
B(t)\le\varepsilon
\quad\Longrightarrow\quad
y_k(t)\ge m,
\]
and, after decreasing $\varepsilon$,
\[
y_i(t),y_j(t)\le\frac m2.
\]
Differentiating $B$, the mutual interaction
between vortices $i$ and $j$ cancels exactly:
\[
s\Gamma_i\Gamma_j(x_i-x_j)
+s\Gamma_j\Gamma_i(x_j-x_i)=0,
\]
only their interactions with vortex $k$ remain, we obtain
\[
\begin{aligned}
B'(t)
=\frac{2s\Gamma_k}{\pi}
\sum_{\ell\in\{i,j\}}
\Gamma_\ell
\frac{(x_\ell-x_k)y_\ell y_k}
{|z_\ell-z_k|^2|z_\ell-\overline z_k|^2}.
\end{aligned}
\]
The same separation estimates as in
Lemma~\ref{lem:isolated} give
\[
|B'(t)|\le CB(t)
\quad\text{whenever }B(t)\le\varepsilon.
\]
By Lemma~\ref{lem:barrier}(i),
\[
B(t)\ge\min\{B(0),\varepsilon\}e^{-Ct}.
\]
This contradicts $B(t_n)\to0$.

Therefore two vortices of the same circulation sign cannot simultaneously approach the boundary either.

Combining the two cases proves that for every pair \(i\ne j\),
\[
\inf_{0\le t<T}\max\{y_i(t),y_j(t)\}>0.
\]
Fixing $i$, we can therefore choose
$\varepsilon,m>0$ such that
\[
y_i(t)\le\varepsilon
\quad\Longrightarrow\quad
y_j(t)\ge m
\quad(j\ne i).
\]
Lemma~\ref{lem:isolated} implies
\[
\inf_{0\le t<T}y_i(t)>0.
\]
Since $i$ is arbitrary, therefore an individual vortex cannot reach the boundary, completing the proof.
\end{proof}

\subsection{Zero vertical impulse: common estimates}

We now assume
\[
P(0)=0.
\]
Because every $y_i>0$ and the circulations
are nonzero, the circulations must have mixed
signs. After relabelling, there exist
$a,b,c>0$ and $s\in\{-1,1\}$ such that
\begin{equation}
\label{eq:circs}
\Gamma_1=sa,\qquad
\Gamma_2=sb,\qquad
\Gamma_3=-sc.
\end{equation}
Set
\[
A=a+b,\qquad \alpha=\frac Ac.
\]
Define
\begin{equation}
\label{eq:center}
Z=X+iY=\frac{az_1+bz_2}{A},
\qquad
w=z_1-z_2,\qquad r=|w|.
\end{equation}
Then
\begin{equation}
\label{eq:inverse}
z_1=Z+\frac bA w,\qquad
z_2=Z-\frac aA w.
\end{equation}
Conservation of the zero vertical impulse gives
\begin{equation}
\label{eq:heightrelation}
ay_1+by_2=cy_3,
\qquad
y_3=\alpha Y.
\end{equation}
The following concentration estimate is common
to both zero-impulse cases.

\begin{lemma}[(Energy concentration of the like-signed pair)]
\label{lem:concentration}
Under \eqref{eq:circs}--\eqref{eq:heightrelation},
there exist constants $C>0$ and
$\varepsilon_0>0$, depending only on the
circulations and the initial energy, such that
\[
r\le CY^\kappa,\qquad
\kappa=\frac{A^2+c^2}{2ab}>2,
\]
whenever $0<Y\le\varepsilon_0$.

In particular,
\[
\frac rY\to0
\quad\text{uniformly on the fixed energy level
as }Y\downarrow0,
\]
and
\[
y_1=Y+O(r),\qquad y_2=Y+O(r).
\]
\end{lemma}

\begin{proof}
Introduce
\[
L_{ij}
=\log\frac{S_{ij}^2}{R_{ij}^2}
=\log\left(
1+\frac{4y_iy_j}{R_{ij}^2}
\right)\ge0.
\]
The conserved energy becomes
\begin{equation}
\label{eq:energy-mixed}
E_0=
\sum_{i=1}^3\Gamma_i^2\log(2y_i)
+abL_{12}-acL_{13}-bcL_{23}.
\end{equation}
By \eqref{eq:heightrelation},
\[
y_1\le\frac AaY,\qquad
y_2\le\frac AbY,\qquad
y_3=\alpha Y.
\]
Let
$\sigma=a^2+b^2+c^2.$
Then,
\[
\sum_{i=1}^3\Gamma_i^2\log(2y_i)
\le\sigma\log Y+C_1.
\]
Since $L_{13},L_{23}\ge0$,
\eqref{eq:energy-mixed} implies
\[
abL_{12}
\ge E_0-\sigma\log Y-C_1.
\]
Thus
\[
L_{12}\ge
\frac{\sigma}{ab}\log\frac1Y-C_2.
\]
For sufficiently small $Y$,
\[
e^{L_{12}}-1
\ge C_3Y^{-\sigma/(ab)}.
\]
Using
$e^{L_{12}}-1=\frac{4y_1y_2}{r^2},$
we obtain
\[
r^2
\le CY^{2+\sigma/(ab)}.
\]
But
$
2+\frac{\sigma}{ab}
=\frac{A^2+c^2}{ab}.
$
Consequently,
\[
r\le CY^{(A^2+c^2)/(2ab)}.
\]
Since
$A^2=(a+b)^2\ge4ab,$
we have
\[
\kappa\ge2+\frac{c^2}{2ab}>2.
\]
Finally, \eqref{eq:inverse} implies
\[
|y_1-Y|\le\frac bA r,
\qquad
|y_2-Y|\le\frac aA r.
\]
The remaining assertions follow.
\end{proof}

We also record a monotonicity identity.

\begin{lemma}[(Horizontal monotonicity)]
\label{lem:horizontal}
Define
\[
Q(t)=\sum_{i=1}^3\Gamma_i x_i(t).
\]
Then
\begin{equation}
\label{eq:Qmonotone}
\dot Q(t)=\frac1{2\pi}\int_0^\infty
\left|
\sum_{j=1}^3
\Gamma_j e^{-ky_j(t)}e^{ikx_j(t)}
\right|^2\,dk>0.
\end{equation}
\end{lemma}

\begin{proof}
Multiplying the horizontal equation in
\eqref{eq:system} by $\Gamma_i$ and summing
over $i$, the direct interactions cancel
pairwise. Thus
\[
\dot Q=
\sum_{i=1}^3\frac{\Gamma_i^2}{4\pi y_i}
+\frac1\pi\sum_{i<j}
\Gamma_i\Gamma_j
\frac{y_i+y_j}{S_{ij}^2}.
\]
On the other hand,
\[
\int_0^\infty e^{-2ky_i}\,dk=\frac1{2y_i},
\]
and
\[
\int_0^\infty
e^{-k(y_i+y_j)}
\cos(k(x_i-x_j))\,dk
=\frac{y_i+y_j}{S_{ij}^2}.
\]
Expanding the squared modulus in
\eqref{eq:Qmonotone} proves the identity
and nonnegativity.

To prove strict positivity, suppose
$\dot Q(t)=0$ at some time. Then
\[
\sum_{j=1}^3\Gamma_j e^{ikz_j(t)}=0
\qquad\text{for every }k\ge0.
\]
Differentiating at $k=0$ for orders
$m=0,1,2$ gives
\[
\sum_{j=1}^3\Gamma_j z_j(t)^m=0,
\qquad m=0,1,2.
\]
Since the three positions are pairwise
distinct, the associated Vandermonde matrix
is invertible. Hence all $\Gamma_j=0$,
contrary to the assumptions.

Therefore $\dot Q(t)>0$.
\end{proof}

\subsection{Zero vertical impulse: zero total circulation}

We first consider the neutral configuration
\[
\Gamma_1+\Gamma_2+\Gamma_3=0.
\]
Under \eqref{eq:circs}, this is equivalent
to $c=A$.
\begin{proposition}[(Neutral three-vortex system)]
\label{prop:neutral}
Suppose
\[
P(0)=0,\qquad
\Gamma_1+\Gamma_2+\Gamma_3=0.
\]
Then, for every finite interval of classical
collision-free existence $[0,T)$,
\[
\inf_{0\le t<T}y_i(t)>0,
\qquad i=1,2,3.
\]
\end{proposition}

\begin{proof}
Under the assumptions,
$c=A,\, y_3=Y.$
Define the horizontal separation
\[
D=X-x_3,
\qquad
W=sD.
\]
We first show that, whenever $Y$ is small,
the separation $|D|$ is much larger than
the internal pair separation $r$.

\medskip
\noindent
\textit{Step 1. Uniform separation.}

We claim that, for every $K>0$, there
exists $\varepsilon_K>0$ such that
\begin{equation}
\label{eq:neutral-separation}
Y\le\varepsilon_K
\quad\Longrightarrow\quad
|D|>Kr.
\end{equation}
Suppose otherwise. Then there exists a
sequence of configurations on the same
energy and impulse level such that
\[
Y_n\to0,\qquad |D_n|\le Kr_n.
\]
By Lemma~\ref{lem:concentration},
$r_n/Y_n\to0$. Moreover,
\[
y_1/Y\to1,\qquad y_2/Y\to1,
\qquad y_3=Y.
\]
From \eqref{eq:inverse},
\[
|z_1-z_3|,\ |z_2-z_3|
\le(K+1)r.
\]
Consequently,
\[
L_{12}=2\log(Y/r)+O(1),
\]
and
\[
L_{13},L_{23}
\ge2\log(Y/r)-C_K.
\]
Also,
\[
\sum_{i=1}^3\Gamma_i^2\log(2y_i)
=\sigma\log Y+O(1),
\]
where
$\sigma=a^2+b^2+A^2.$
Substitution into \eqref{eq:energy-mixed}
gives
\[
E_0\le
\sigma\log Y
+2(ab-A^2)\log(Y/r)+C_K'.
\]
Since
$
ab-A^2<0,\,
Y\to0,\, Y/r\to\infty,
$
the right-hand side tends to $-\infty$,
contradicting the fixed finite energy $E_0$.
This proves \eqref{eq:neutral-separation}.

In particular, there exists $\varepsilon>0$
such that
\begin{equation}
\label{eq:Ddominates}
Y\le\varepsilon
\quad\Longrightarrow\quad
r<\frac12|D|.
\end{equation}

\medskip
\noindent
\textit{Step 2. Monotonicity of $W$.}

Since
\[
Q=sa x_1+sb x_2-sA x_3
=sA(X-x_3)=AW,
\]
Lemma~\ref{lem:horizontal} yields
\begin{equation}
\label{eq:Wpositive}
\dot W=\frac{\dot Q}{A}>0.
\end{equation}
Thus $W$ is strictly increasing.

\medskip
\noindent
\textit{Step 3. Evolution of the pair height.}

Differentiating
$
Y=\frac{ay_1+by_2}{A},
$
the internal interaction between vortices
1 and 2 cancels. The vertical equation
in \eqref{eq:system}, together with
$\Gamma_3=-sA$ and $y_3=Y$, gives
\begin{equation}
\label{eq:Yneutral}
\dot Y=
-\frac{2sY}{\pi}
\sum_{i=1}^2
\ell_i
\frac{(x_i-x_3)y_i}{R_{i3}^2S_{i3}^2},
\end{equation}
where,
$
\ell_1=a,\, \ell_2=b.
$
By \eqref{eq:inverse},
\[
x_1-x_3=D+\frac bA\operatorname{Re}w,
\]
\[
x_2-x_3=D-\frac aA\operatorname{Re}w.
\]
If $Y\le\varepsilon$, then
\eqref{eq:Ddominates} implies
\[
\operatorname{sgn}(x_i-x_3)
=\operatorname{sgn}D,
\]
and
\begin{equation}
\label{eq:horizontal-distance}
|x_i-x_3|\ge\frac12|D|
=\frac12|W|,
\qquad i=1,2.
\end{equation}
It follows from \eqref{eq:Yneutral} that
\begin{equation}
\label{eq:Ysign}
\operatorname{sgn}\dot Y
=-\operatorname{sgn}W
\quad\text{whenever }Y\le\varepsilon.
\end{equation}

\medskip
\noindent
\textit{Step 4. A positive lower bound for $Y$.}

We distinguish two possibilities.

\smallskip
\noindent
\emph{Case I: $W(t)<0$ for all $0\le t<T$.}

Equation \eqref{eq:Ysign} gives
\[
\dot Y(t)>0
\quad\text{whenever }Y(t)\le\varepsilon.
\]
Applying Lemma~\ref{lem:barrier}(i)
with $C=0$ yields
\[
Y(t)\ge\min\{Y(0),\varepsilon\}>0.
\]

\smallskip
\noindent
\emph{Case II: $W(t_0)>0$ for some $t_0<T$.}

Since $W$ is strictly increasing,
\[
W(t)\ge W(t_0)=:d_0>0
\qquad(t_0\le t<T).
\]
Whenever $Y(t)\le\varepsilon$,
\eqref{eq:horizontal-distance} implies
$
|x_i-x_3|\ge\frac{d_0}{2},
\; i=1,2.
$
Using
$
R_{i3},S_{i3}\ge|x_i-x_3|,
$
we obtain
\[
\frac{|x_i-x_3|}
{R_{i3}^2S_{i3}^2}
\le\frac8{d_0^3}.
\]
Equation \eqref{eq:Yneutral} therefore gives
\[
|\dot Y|
\le\frac{16Y}{\pi d_0^3}
(ay_1+by_2)
=\frac{16A}{\pi d_0^3}Y^2.
\]
In particular,
\[
\dot Y\ge-C_0Y^2,
\qquad
C_0=\frac{16A}{\pi d_0^3},
\]
whenever $Y\le\varepsilon$.

Applying Lemma~\ref{lem:barrier}(ii)
on $[t_0,T)$ gives
\[
\frac1{Y(t)}
\le
\max\left\{\frac1{Y(t_0)},\frac1\varepsilon\right\}
+C_0(t-t_0).
\]
Because $T<\infty$, $Y$ is bounded
away from zero on $[t_0,T)$.
It is also bounded away from zero on
$[0,t_0]$ by continuity.

These two cases exhaust all possibilities:
if a strictly increasing $W$ is not
negative throughout $[0,T)$, then it
is positive at some time before $T$.

We have proved
\begin{equation}
\label{eq:Yneutral-lower}
\inf_{0\le t<T}Y(t)=:\eta>0.
\end{equation}

\medskip
\noindent
\textit{Step 5. Individual vortex heights.}

Because $y_3=Y$, we immediately have
\[
y_3(t)\ge\eta.
\]
Moreover,
\[
ay_1+by_2=AY\ge A\eta.
\]
Whenever $y_1$ is sufficiently small,
this identity forces $y_2$ to remain
bounded away from zero. Since $y_3$
is also bounded below, the hypotheses
of Lemma~\ref{lem:isolated} hold for
vortex 1. Hence
\[
\inf_{0\le t<T}y_1(t)>0.
\]
Interchanging indices 1 and 2 gives
\[
\inf_{0\le t<T}y_2(t)>0.
\]
This proves the proposition.
\end{proof}

\subsection{Nonzero total circulation}

We now consider the remaining zero-impulse
configuration
$
\Gamma_1+\Gamma_2+\Gamma_3\ne0.
$
In the notation \eqref{eq:circs}, this
is equivalent to $A\ne c$.

Unlike the neutral case, the horizontal
separation $s(X-x_3)$ need not be monotone.
Instead, we introduce its ratio to the
pair height $Y$ and prove that this
dimensionless variable increases whenever
$Y$ is sufficiently small.

\begin{proposition}[(Non-neutral three-vortex system)]
\label{prop:nonneutral}
Suppose
\[
P(0)=0,\qquad
\Gamma_1+\Gamma_2+\Gamma_3\ne0.
\]
Then, for every finite interval of classical
collision-free existence $[0,T)$,
\[
\inf_{0\le t<T}y_i(t)>0,
\qquad i=1,2,3.
\]
\end{proposition}

\begin{proof}
Write
$
\Gamma_1=sa,\;
\Gamma_2=sb,\;
\Gamma_3=-sc,
$
where $a,b,c>0$, $s\in\{-1,1\}$,
and $A=a+b\ne c$.
Recall that
\[
y_3=\alpha Y,\qquad
\alpha=\frac Ac\ne1.
\]

\medskip
\noindent
\textit{Step 1. Dimensionless coordinates.}

Define
\[
W=s(X-x_3),\qquad h=\frac WY,
\]
and
\[
q=\frac{s\operatorname{Re}w}{Y},
\qquad
\nu=\frac{\operatorname{Im}w}{Y},
\qquad
\rho=\sqrt{q^2+\nu^2}=\frac rY.
\]
Let
\[
\lambda_1=\frac bA,\qquad
\lambda_2=-\frac aA,
\]
and put
\[
u_i=1+\lambda_i\nu,\qquad
p_i=h+\lambda_i q,
\qquad i=1,2.
\]
By \eqref{eq:inverse},
\[
y_i=Yu_i,\qquad
s(x_i-x_3)=Yp_i.
\]
Define
\[
\mathcal R_i=p_i^2+(u_i-\alpha)^2,
\qquad
\mathcal S_i=p_i^2+(u_i+\alpha)^2.
\]
Then
\[
R_{i3}^2=Y^2\mathcal R_i,
\qquad
S_{i3}^2=Y^2\mathcal S_i.
\]
Since $\alpha\ne1$, we can choose
$\rho_0>0$ sufficiently small that
\[
\frac34\le u_i\le\frac54,
\qquad
|u_i-\alpha|\ge\frac{|1-\alpha|}{2}
\]
whenever $\rho\le\rho_0$.

Moreover, since
\[
|p_i-h|\le\rho_0,
\]
there exists $m>0$, depending only on
$a,b,c$, such that
\begin{equation}
\label{eq:denominator-bounds}
\mathcal R_i,\mathcal S_i
\ge m(1+h^2)
\end{equation}
for every $h\in\mathbb R$ and
$\rho\le\rho_0$.

Indeed, when $|h|\ge2\rho_0$,
$|p_i|\ge|h|/2$; when
$|h|<2\rho_0$, the positive lower
bounds on $|u_i-\alpha|$ and
$u_i+\alpha$ apply.

\medskip
\noindent
\textit{Step 2. Exact scaled equations.}

Differentiating $Y=(ay_1+by_2)/A$,
the interactions between vortices 1
and 2 cancel. The vertical equations
give
\begin{equation}
\label{eq:scaled-Y}
Y\dot Y=\mathcal V(h,q,\nu),
\end{equation}
where
\begin{equation}
\label{eq:V-def}
\mathcal V(h,q,\nu)
=-\frac2\pi
\sum_{i=1}^2
\ell_i
\frac{p_i u_i}{\mathcal R_i\mathcal S_i},
\end{equation}
with $\ell_1=a,\ell_2=b$.

We next derive the horizontal equation
for $W=s(X-x_3)$.

Using the horizontal equations in
\eqref{eq:system}, the singular direct
interaction between vortices 1 and 2
cancels in the weighted velocity
\[
\dot X=\frac{a\dot x_1+b\dot x_2}{A}.
\]
A direct calculation gives
\begin{align}
\dot W={}&
\frac1{4\pi A}
\left(\frac{a^2}{y_1}+\frac{b^2}{y_2}\right)
+\frac{ab}{\pi A}
\frac{y_1+y_2}{S_{12}^2}
+\frac{c}{4\pi y_3}
\nonumber\\
&-\frac1{2\pi}\sum_{i=1}^2\ell_i
\left[
\left(1+\frac cA\right)
\frac{y_i+y_3}{S_{i3}^2}
+
\left(1-\frac cA\right)
\frac{y_i-y_3}{R_{i3}^2}
\right].
\label{eq:W-exact}
\end{align}
Substituting the dimensionless variables,
we obtain
\begin{equation}
\label{eq:scaled-W}
Y\dot W=\mathcal G(h,q,\nu),
\end{equation}
where
\begin{align}
\mathcal G
={}&
\frac1{4\pi A}
\left(\frac{a^2}{u_1}+\frac{b^2}{u_2}\right)
+\frac{ab}{\pi A}
\frac{u_1+u_2}
{q^2+(u_1+u_2)^2}
+\frac{c}{4\pi\alpha}
\nonumber\\
&-\frac1{2\pi}\sum_{i=1}^2\ell_i
\left[
\left(1+\frac cA\right)
\frac{u_i+\alpha}{\mathcal S_i}
+
\left(1-\frac cA\right)
\frac{u_i-\alpha}{\mathcal R_i}
\right].
\label{eq:G-def}
\end{align}
Both \eqref{eq:scaled-Y} and
\eqref{eq:scaled-W} are exact identities.

\medskip
\noindent
\textit{Step 3. Uniform positivity in the
tight-pair limit.}

At $q=\nu=0$, we have
\[
u_1=u_2=1,\qquad p_1=p_2=h.
\]
Formula \eqref{eq:G-def} reduces to
\begin{align}
\mathcal G(h,0,0)
={}&
\frac{A^2+c^2}{4\pi A}
-\frac{(A+c)(1+\alpha)}
{2\pi[h^2+(1+\alpha)^2]}
\nonumber\\
&+
\frac{(A-c)^2}
{2\pi c[h^2+(1-\alpha)^2]}.
\label{eq:G-limit}
\end{align}
Since the last term is nonnegative and
\[
\frac{(A+c)(1+\alpha)}
{2\pi[h^2+(1+\alpha)^2]}
\le
\frac{A+c}{2\pi(1+\alpha)}
=\frac c{2\pi},
\]
we conclude that
\begin{equation}
\label{eq:G-positive}
\mathcal G(h,0,0)\ge
\frac{(A-c)^2}{4\pi A}
=:g_0>0
\end{equation}
for every $h\in\mathbb R$.

We next prove that this lower bound
persists for sufficiently small $\rho$.

By \eqref{eq:denominator-bounds}, all
the rational functions in
\eqref{eq:G-def} are smooth in
$(q,\nu)$ and have uniformly bounded
first derivatives on
$\{(h,q,\nu):h\in\mathbb R,\ \rho\le\rho_0\}.$


Consider
\[
F_\pm(p,u)
=\frac{u\pm\alpha}
{p^2+(u\pm\alpha)^2}.
\]
Its derivatives are
\[
\partial_pF_\pm
=-\frac{2p(u\pm\alpha)}
{[p^2+(u\pm\alpha)^2]^2},
\]
and
\[
\partial_uF_\pm
=\frac{p^2-(u\pm\alpha)^2}
{[p^2+(u\pm\alpha)^2]^2}.
\]
Using
$
2|p(u\pm\alpha)|
\le p^2+(u\pm\alpha)^2
$
and
$
|p^2-(u\pm\alpha)^2|
\le p^2+(u\pm\alpha)^2,
$
we obtain
\[
|\partial_pF_\pm|,\,
|\partial_uF_\pm|
\le
\frac1{p^2+(u\pm\alpha)^2}.
\]
By \eqref{eq:denominator-bounds},
\[
|\partial_pF_\pm(p_i,u_i)|,\,
|\partial_uF_\pm(p_i,u_i)|
\le\frac1{m(1+h^2)}.
\]
Since
\[
\partial_qp_i=\partial_\nu u_i=\lambda_i,
\qquad
\partial_\nu p_i=\partial_q u_i=0,
\]
the chain rule yields
\[
|\partial_qF_\pm(p_i,u_i)|,\,
|\partial_\nu F_\pm(p_i,u_i)|
\le\frac1m.
\]

The remaining terms in \eqref{eq:G-def}
also have uniformly bounded derivatives.
Indeed, $u_i\ge3/4$ implies
\[
\left|\partial_\nu\frac1{u_i}\right|
\le\frac{16}{9},
\]
while, writing $U=u_1+u_2\ge3/2$,
\[
\left|
\partial_q\frac{U}{q^2+U^2}
\right|
\le\frac49,
\qquad
\left|
\partial_\nu\frac{U}{q^2+U^2}
\right|
\le\frac49.
\]
Consequently, there exists
$C=C(a,b,c)>0$ such that
\[
|\partial_q\mathcal G|
+|\partial_\nu\mathcal G|
\le C
\]
uniformly for $h\in\mathbb R$
and $\rho\le\rho_0$.

Since the disk $\rho\le\rho_0$ is
convex, the mean-value theorem gives
\begin{equation}
\label{eq:G-lipschitz}
|\mathcal G(h,q,\nu)-\mathcal G(h,0,0)|
\le C\rho,
\end{equation}
where the constant is independent of $h$.

Similarly, \eqref{eq:V-def} and
\eqref{eq:denominator-bounds} give
\begin{equation}
\label{eq:V-bound}
|\mathcal V(h,q,\nu)|
\le C\frac{|h|+\rho}{(1+h^2)^2}.
\end{equation}

Decreasing $\rho_0$ if necessary,
\eqref{eq:G-positive} and
\eqref{eq:G-lipschitz} imply
\begin{equation}
\label{eq:G-half}
\mathcal G(h,q,\nu)\ge\frac{g_0}{2}
\end{equation}
whenever $\rho\le\rho_0$.

We claim that
\begin{equation}
\label{eq:GminusV}
\mathcal G-h\mathcal V\ge\frac{g_0}{4}.
\end{equation}

If $|h|>\rho$, then
$
p_i=h+\lambda_i q
$
has the same sign as $h$, because
$|\lambda_iq|\le\rho$.
Since $u_i>0$, \eqref{eq:V-def}
implies
\[
h\mathcal V\le0.
\]
Thus \eqref{eq:GminusV} follows from
\eqref{eq:G-half}.

If $|h|\le\rho$, then
\eqref{eq:V-bound} gives
\[
|h\mathcal V|\le C\rho^2.
\]
We choose $\rho_0$ still smaller so
that $C\rho_0^2\le g_0/4$.
Combining this with \eqref{eq:G-half}
again proves \eqref{eq:GminusV}.

\medskip
\noindent
\textit{Step 4. Monotonicity of the scaled
horizontal separation.}

Recall that
\[
h=\frac WY.
\]
By \eqref{eq:scaled-Y} and
\eqref{eq:scaled-W},
\begin{equation}
\label{eq:h-evolution}
\dot h
=\frac{Y\dot W-W\dot Y}{Y^2}
=\frac{\mathcal G-h\mathcal V}{Y^2}.
\end{equation}
Lemma~\ref{lem:concentration} implies
\[
\rho=\frac rY\le CY^{\kappa-1}
\longrightarrow0
\quad\text{uniformly as }Y\downarrow0.
\]
Hence there exists $\varepsilon>0$
such that
\[
Y\le\varepsilon
\quad\Longrightarrow\quad
\rho\le\rho_0.
\]
By \eqref{eq:GminusV} and
\eqref{eq:h-evolution},
\begin{equation}
\label{eq:h-monotone}
\dot h\ge\frac{g_0}{4Y^2}>0
\quad\text{whenever }Y\le\varepsilon.
\end{equation}
Thus $h$ is strictly increasing on every
time interval during which $Y\le\varepsilon$.

\medskip
\noindent
\textit{Step 5. Integrable variation of
$\log Y$.}

On an interval where $Y\le\varepsilon$,
we may use $h$ as the independent variable.

By \eqref{eq:scaled-Y} and
\eqref{eq:h-evolution},
\[
\frac{d\log Y}{dh}
=
\frac{\mathcal V}{\mathcal G-h\mathcal V}.
\]
Equations \eqref{eq:V-bound} and
\eqref{eq:GminusV} yield
\begin{equation}
\label{eq:logY-bound}
\left|\frac{d\log Y}{dh}\right|
\le
C\frac{1+|h|}{(1+h^2)^2}.
\end{equation}
The right-hand side is integrable over
the whole real line:
\[
\int_{-\infty}^{\infty}
\frac{1+|h|}{(1+h^2)^2}\,dh
=\frac\pi2+1.
\]
Set
\[
M=C\left(\frac\pi2+1\right).
\]
Let $[t_1,t_2]$ be any interval on which
$Y(t)\le\varepsilon$. Since $h$ is strictly
increasing there, integrating
\eqref{eq:logY-bound} gives
\[
\begin{aligned}
|\log Y(t_2)-\log Y(t_1)|
&\le C\int_{h(t_1)}^{h(t_2)}
\frac{1+|h|}{(1+h^2)^2}\,dh\\
&\le M.
\end{aligned}
\]
Consequently,
\begin{equation}
\label{eq:Y-excursion}
Y(t_2)\ge e^{-M}Y(t_1).
\end{equation}
The constant $M$ is independent of
the interval and of the values attained
by $h$.

We now turn this into a uniform lower bound.
Fix $t<T$ with $Y(t)<\varepsilon$.
If $Y$ previously reached $\varepsilon$,
let $\tau\le t$ be its last crossing time.
Then
\[
Y(\tau)=\varepsilon,
\qquad
Y(s)<\varepsilon
\quad(\tau<s\le t).
\]
Equation \eqref{eq:Y-excursion} implies
\[
Y(t)\ge\varepsilon e^{-M}.
\]
If no previous crossing exists, then
$Y(0)<\varepsilon$, and
\eqref{eq:Y-excursion} applied from time
zero gives
\[
Y(t)\ge Y(0)e^{-M}.
\]
When $Y(t)\ge\varepsilon$, a lower bound
is immediate. Therefore
\begin{equation}
\label{eq:Y-nonneutral-lower}
Y(t)\ge e^{-M}\min\{Y(0),\varepsilon\}
=:\eta>0
\end{equation}
for every $0\le t<T$.

\medskip
\noindent
\textit{Step 6. Individual vortex heights.}

Since $y_3=\alpha Y$,
\[
y_3(t)\ge\alpha\eta>0.
\]
Moreover,
\[
ay_1+by_2=AY\ge A\eta.
\]
Whenever $y_1$ is sufficiently small,
the latter relation forces $y_2$ to
remain bounded below by a positive
constant. Together with the lower
bound for $y_3$, this allows us to
apply Lemma~\ref{lem:isolated} to $y_1$.

Hence
\[
\inf_{0\le t<T}y_1(t)>0.
\]
The same argument gives
\[
\inf_{0\le t<T}y_2(t)>0.
\]
Combined with the bound for $y_3$,
this completes the proof.
\end{proof}

\subsection{Proof of the main theorem}

\begin{proof}[Proof of Theorem~\ref{thm:main_bd_noncolli}]
Let $\Gamma_1,\Gamma_2,\Gamma_3$ be arbitrary
nonzero real circulations.

By conservation of vertical impulse,
either $P(0)\ne0$ or $P(0)=0$.

If $P(0)\ne0$, Proposition~\ref{prop:Pnonzero}
establishes the required positive lower
bounds for all three vortex heights.

Suppose $P(0)=0$. If
$
\Gamma_1+\Gamma_2+\Gamma_3=0,
$
the assertion follows from
Proposition~\ref{prop:neutral}.

If instead
$
\Gamma_1+\Gamma_2+\Gamma_3\ne0,
$
it follows from
Proposition~\ref{prop:nonneutral}.

These cases exhaust all possible
circulation configurations.

Therefore, for every finite interval
of classical collision-free existence,
\[
\inf_{0\le t<T}y_i(t)>0,
\qquad i=1,2,3.
\]
This proves the theorem.
\end{proof}

\begin{remark}[(Boundary versus interior collapse)]
Theorem~\ref{thm:main_bd_noncolli} concerns only
collisions with the boundary. It does
not exclude finite-time collisions
between vortices in the interior
of $\mathbb H^+$.

Furthermore, the proof makes no
assumption that the heights or vortex
positions possess limits as $t\uparrow T$.
All conclusions hold in the standard
$\liminf$ sense.
\end{remark}


\medskip

\section{Rigorous Existence of Finite-Time Interior Collapses in the Half-Plane}
\label{sec:rigorous_existence_collapse}

In this section, we construct a rigorous existence proof for finite-time interior collapses of $N=3$ point vortices in the upper half-plane $\mathbb H^+$. 


\subsection{Formal Derivation and Unperturbed Profile Equations}

Consider three point vortices at positions $z_k(t) = x_k(t) + i y_k(t) \in \mathbb{C}^+,\,k=1,\,2,\,3,$ with circulations $\Gamma_k \in \mathbb{R} \setminus \{0\}$ satisfying the Gr\"{o}bli self-similar collapse condition:
$$
\Gamma_1 \Gamma_2 + \Gamma_2 \Gamma_3 + \Gamma_3 \Gamma_1 = 0.
$$
Because $\sum_{j=1}^3 \Gamma_j \neq 0$ for non-trivial choices (e.g., $\boldsymbol{\Gamma} = (2, 2, -1)$), we define the instantaneous center of vorticity:
\begin{equation}
C(t) = \frac{1}{\sum_{j=1}^3 \Gamma_j} \sum_{j=1}^3 \Gamma_j z_j(t).
\end{equation}
By Newton's Third Law, the singular Biot--Savart interactions conserve the total circulation vector, so $C(t)$ evolves exclusively under the smooth regular boundary field:
\begin{equation}
\label{eq:center_evolution}
\frac{d\bar{C}}{dt} = \frac{1}{\sum_{j=1}^3 \Gamma_j} \sum_{j=1}^3 \Gamma_j \overline{V_{\text{reg}}^{(j)}(z_j, \mathbf{z})}, \quad \text{where } V_{\text{reg}}^{(j)}(z_j, \mathbf{z}) = \frac{1}{2\pi i} \sum_{m=1}^3 \frac{\Gamma_m}{\bar{z}_j - z_m}.
\end{equation}

In terms of the relative coordinates $\zeta_k(t) = z_k(t) - C(t)$, the governing equations take the split form:
\begin{equation}
\label{eq:relative_system}
\frac{d\bar{\zeta}_k}{dt} = \overline{V_{\text{sing}}^{(k)}(\boldsymbol{\zeta})} + \left( \overline{V_{\text{reg}}^{(k)}(z_k, \mathbf{z})} - \frac{d\bar{C}}{dt} \right), \quad \text{with } V_{\text{sing}}^{(k)}(\boldsymbol{\zeta}) = -\frac{1}{2\pi i} \sum_{\substack{j=1 \\ j \neq k}}^3 \frac{\Gamma_j}{\bar{\zeta}_k - \bar{\zeta}_j}.
\end{equation}

\begin{remark}[($V_{\text{sing}}^{(k)}(\boldsymbol{\zeta})$ vs $V_{\text{sing}}^{(k)}(\mathbf{z})$)]

Because $\zeta_k = z_k - C$, the relative displacement vectors satisfy:
$$\zeta_k - \zeta_j = (z_k - C) - (z_j - C) = z_k - z_j.$$
Therefore:
$V_{\text{sing}}^{(k)}(\mathbf{z}) \equiv V_{\text{sing}}^{(k)}(\boldsymbol{\zeta}).$
Writing $V_{\text{sing}}^{(k)}(\boldsymbol{\zeta})$ emphasizes that the singular Biot--Savart interactions depend purely on internal relative distances, making them invariant under spatial translations.

Unlike $V_{\text{sing}}$, the regular image field $V_{\text{reg}}^{(k)}(z_k, \mathbf{z})$ depends on absolute positions $z_k = C + \zeta_k$, because images are anchored to the boundary $\partial \mathbb H^+$. Expressing $V_{\text{sing}}$ in $\boldsymbol{\zeta}$ highlights this physical distinction.
\end{remark}

Neglecting the boundary correction in \eqref{eq:relative_system} gives the unperturbed full-plane system. Substituting the self-similar collapsing form
$$
\zeta_k^{(0)}(t)=\sqrt{T-t}\,e^{i\Omega\ln(T-t)}W_k
$$
into this system yields the nonlinear algebraic equation for the profile $\mathbf{W} = (W_1, W_2, W_3)^T \in \mathbb{C}^3$:
\begin{equation}
\label{eq:profile_algebraic}
\left( \frac{1}{2} - i\Omega \right) \bar{W}_k = -\frac{1}{2\pi i} \sum_{\substack{j=1 \\ j \neq k}}^3 \frac{\Gamma_j}{W_k - W_j}.
\end{equation}


\begin{proof}[Derivation of Equation \eqref{eq:profile_algebraic}]
We consider the unperturbed full-plane $N$-vortex system expressed in relative coordinates $\zeta_k(t) = z_k(t) - C(t)$, given by
\begin{equation}
\label{eq:proof_relative_system}
\frac{d\bar{\zeta}_k}{dt} = \frac{1}{2\pi i} \sum_{\substack{j=1 \\ j \neq k}}^3 \frac{\Gamma_j}{\zeta_k - \zeta_j}.
\end{equation}

Let $\tau = T - t > 0$ denote the time remaining until the finite-time collapse at $t = T$, then $\frac{d}{dt} = -\frac{d}{d\tau}.$ 
We insert the self-similar collapsing form into relative coordinates:
\begin{equation}
\label{eq:proof_ansatz}
\zeta_k(\tau) = \sqrt{\tau} \, e^{i\Omega \ln\tau} W_k \implies \bar{\zeta}_k(\tau) = \sqrt{\tau} \, e^{-i\Omega \ln\tau} \bar{W}_k,
\end{equation}
where $W_k \in \mathbb{C}$ represents the scale-invariant shape coordinate and $\Omega \in \mathbb{R}$ is the asymptotic rotation rate.

Tthe left-hand side of \eqref{eq:proof_relative_system} is given by
\begin{equation}
\label{eq:proof_lhs}
\frac{d\bar{\zeta}_k}{dt} = -\frac{d\bar{\zeta}_k}{d\tau} = -\frac{1}{\sqrt{\tau}} e^{-i\Omega \ln\tau} \left( \frac{1}{2} - i\Omega \right) \bar{W}_k.
\end{equation}

Next, substituting \eqref{eq:proof_ansatz} into the right-hand side of \eqref{eq:proof_relative_system} yields:
\begin{align}
\label{eq:proof_rhs}
\frac{1}{2\pi i} \sum_{j \neq k}^3 \frac{\Gamma_j}{\zeta_k - \zeta_j} &= \frac{1}{2\pi i} \sum_{j \neq k}^3 \frac{\Gamma_j}{\sqrt{\tau} e^{i\Omega \ln\tau} W_k - \sqrt{\tau} e^{i\Omega \ln\tau} W_j} \nonumber \\
&= \frac{1}{\sqrt{\tau}} e^{-i\Omega \ln\tau} \left( \frac{1}{2\pi i} \sum_{j \neq k}^3 \frac{\Gamma_j}{W_k - W_j} \right).
\end{align}

Equating the expressions for the left-hand side \eqref{eq:proof_lhs} and right-hand side \eqref{eq:proof_rhs}, we obtain:
$$
-\frac{1}{\sqrt{\tau}} e^{-i\Omega \ln\tau} \left( \frac{1}{2} - i\Omega \right) \bar{W}_k = \frac{1}{\sqrt{\tau}} e^{-i\Omega \ln\tau} \left( \frac{1}{2\pi i} \sum_{j \neq k}^3 \frac{\Gamma_j}{W_k - W_j} \right).
$$

 Canceling the common factor $\frac{1}{\sqrt{\tau}} e^{-i\Omega \ln\tau} \neq 0$ from both sides yields the algebraic profile equation,
which matches equation \eqref{eq:profile_algebraic} identically.
\end{proof}

Equilateral configurations reduce \eqref{eq:profile_algebraic} to a purely rotational rigid equilibrium.  
Thus, a collapsing profile \(\mathbf W^*\) must have a non-equilateral shape (a scalene triangle or a collinear configuration) that allows the three vortices to contract toward a common point while rotating together. The real and imaginary parts of the profile equation determine the contraction rate and angular velocity, respectively.

\subsection{Similarity Variables and Regularity of the boundary remainder}

To map the singular limit $t\uparrow T$ to an asymptotic steady-state analysis as $s \to \infty$, we introduce similarity variables $(s, \mathbf{Y}(s))$:
$$
s = -\ln(T-t) \implies \frac{ds}{dt} = e^s, \qquad Y_k(s) = e^{s/2} e^{i\Omega s} \zeta_k(t),\,k=1,\,2,\,3.
$$
Differentiating $\bar{Y}_k(s)$ with respect to $s$ and utilizing the homogeneous scale invariance of degree $-1$ of $V_{\text{sing}}^{(k)}$ yields:
\begin{equation}
\label{eq:similarity_ode}
\frac{d\bar{Y}_k}{ds} = \left( \frac{1}{2} - i\Omega \right) \bar{Y}_k + \overline{V_{\text{sing}}^{(k)}(\mathbf{Y})} + R_k(\mathbf{Y}, C, s),
\end{equation}
where the boundary remainder term $R_k$ is defined by:
\begin{equation}\label{eq:similarity_remainder}
R_k(\mathbf{Y}, C, s) = e^{-s/2} e^{-i\Omega s} \left( \overline{V_{\text{reg}}^{(k)}(z_k, \mathbf{z})} - \frac{d\bar{C}}{dt} \right).
\end{equation}

\begin{proof}[Derivation of Similarity ODE \eqref{eq:similarity_ode}]
Taking the complex conjugate of $Y_k(s)$ gives
$$
\bar{Y}_k(s) = e^{(1/2 - i\Omega)s} \bar{\zeta}_k(t) \iff \bar{\zeta}_k(t) = e^{-(1/2 - i\Omega)s} \bar{Y}_k(s).
$$

Differentiating $\bar{Y}_k(s)$ with respect to similarity time $s$ using the product rule yields:
$$
\frac{d\bar{Y}_k}{ds}= \left( \frac{1}{2} - i\Omega \right) e^{(1/2 - i\Omega)s} \bar{\zeta}_k(t) + e^{(1/2 - i\Omega)s} \frac{d\bar{\zeta}_k}{dt} \frac{dt}{ds}.
$$
Substituting $\bar{Y}_k(s) = e^{(1/2 - i\Omega)s} \bar{\zeta}_k(t)$ and $\frac{dt}{ds} = e^{-s}$, the derivative simplifies to:
\begin{equation}
\label{eq:proof_dY_ds_intermediate}
\frac{d\bar{Y}_k}{ds} = \left( \frac{1}{2} - i\Omega \right) \bar{Y}_k + e^{-(1/2 + i\Omega)s} \frac{d\bar{\zeta}_k}{dt}.
\end{equation}

Next, we substitute the full physical system equation \eqref{eq:relative_system} for $\frac{d\bar{\zeta}_k}{dt}$ into \eqref{eq:proof_dY_ds_intermediate}:
\begin{equation}
\label{eq:proof_sub_relative}
\frac{d\bar{Y}_k}{ds} = \left( \frac{1}{2} - i\Omega \right) \bar{Y}_k + e^{-(1/2 + i\Omega)s} \overline{V_{\text{sing}}^{(k)}(\boldsymbol{\zeta})} + e^{-(1/2 + i\Omega)s} \left( \overline{V_{\text{reg}}^{(k)}(z_k, \mathbf{z})} - \frac{d\bar{C}}{dt} \right).
\end{equation}

To see why no prefactor remains before $\overline{V_{\text{sing}}^{(k)}(\mathbf{Y})}$, observe that the Biot--Savart kernel $V_{\text{sing}}^{(k)}$ is homogeneous of degree $-1$. Using $\zeta_m = e^{-(1/2 + i\Omega)s} Y_m$, we rewrite the singular interaction sum in terms of $\mathbf{Y}$:
\begin{align}
\overline{V_{\text{sing}}^{(k)}(\boldsymbol{\zeta})} &= \frac{1}{2\pi i} \sum_{j \neq k}^3 \frac{\Gamma_j}{\zeta_k - \zeta_j} \nonumber 
= \frac{1}{2\pi i} \sum_{j \neq k}^3 \frac{\Gamma_j}{e^{-(1/2 + i\Omega)s} (Y_k - Y_j)} \nonumber \\
&= e^{(1/2 + i\Omega)s} \left( \frac{1}{2\pi i} \sum_{j \neq k}^3 \frac{\Gamma_j}{Y_k - Y_j} \right)
= e^{(1/2 + i\Omega)s} \overline{V_{\text{sing}}^{(k)}(\mathbf{Y})}.\nonumber
\end{align}
Substituting this scaling identity back into equation \eqref{eq:proof_sub_relative} produces an exact cancellation of exponential factors.

Defining the regular boundary remainder as
$$
R_k(\mathbf{Y}, C, s) := e^{-(1/2 + i\Omega)s} \left( \overline{V_{\text{reg}}^{(k)}(z_k, \mathbf{z})} - \frac{d\bar{C}}{dt} \right),
$$
we obtain the final similarity system \eqref{eq:similarity_ode}.
\end{proof}

Next we study the regularity of the boundary remainder $\widetilde R_j(\mathbf Y,C,\eta):=\overline{R_j(\mathbf Y,C,-2\log\eta)}$, where $\eta=e^{-s/2}>0$.
The regularity proposition will be used in Section~\ref{subsec:interior_collapse}.

\begin{proposition}[(Regularity of the boundary remainder)]
\label{prop:boundary-remainder-C1}
Let $\Gamma_1,\Gamma_2,\Gamma_3\in\mathbb R$ satisfy
$\Gamma_{\mathrm{tot}}:=\sum_{j=1}^3\Gamma_j\ne0$,
and let $\Omega\in\mathbb R$. Fix
$\mathbf W^*\in\mathbb C^3$ and $C_*\in\mathbb C$ with
$\operatorname{Im}C_*>0$.
For $\eta=e^{-s/2}>0$, set
\[
q(\eta):=e^{2i\Omega\log\eta},
\qquad
z_j:=C+\eta q(\eta)Y_j,
\qquad
D_{jm}:=z_j-\bar z_m.
\]
Define
\[
U_j:=\overline{V_{\mathrm{reg}}^{(j)}}
=-\frac{1}{2\pi i}\sum_{m=1}^3\frac{\Gamma_m}{D_{jm}},
\qquad
U_c:=\frac{1}{\Gamma_{\mathrm{tot}}}
\sum_{\ell=1}^3\Gamma_\ell U_\ell
=\frac{d\bar C}{dt},
\]
and
\[
R_j(\mathbf Y,C,s)
:=e^{-(1/2+i\Omega)s}(U_j-U_c),
\qquad
\widetilde R_j(\mathbf Y,C,\eta)
:=\overline{R_j(\mathbf Y,C,-2\log\eta)}.
\]

There exist bounded convex neighborhoods $\mathcal U_Y$ of
$\mathbf W^*$ and $\mathcal U_C$ of $C_*$, and constants
$\eta_0\in(0,1)$ and $K>0$, such that, uniformly on
$\mathcal U_Y\times\mathcal U_C\times(0,\eta_0)$,
\begin{align}
\|\widetilde{\mathbf R}(\mathbf Y,C,\eta)\|
&\le K\eta^2\operatorname{diam}(\mathbf Y)
\le 2K\eta^2\|\mathbf Y\|,
\label{eq:boundary-remainder-size}\\
\|D_{\mathbf Y}\widetilde{\mathbf R}\|
+\|D_C\widetilde{\mathbf R}\|
&\le K\eta^2,
\qquad
\|\partial_\eta\widetilde{\mathbf R}\|\le K\eta.
\label{eq:boundary-remainder-derivatives}
\end{align}
Here
\[
\|\mathbf Y\|:=\max_j|Y_j|,
\qquad
\operatorname{diam}(\mathbf Y):=
\max_{j,\ell}|Y_j-Y_\ell|,
\]
and derivatives with respect to complex variables are understood
as real derivatives.

Consequently, setting $\widetilde{\mathbf R}=0$ for $\eta\le0$
gives a $C^1$ extension across $\eta=0$, satisfying
\[
\widetilde{\mathbf R}(\mathbf Y,C,0)=0,
\qquad
D\widetilde{\mathbf R}(\mathbf Y,C,0)=0.
\]
The original remainder satisfies, for
$s\ge s_0:=-2\log\eta_0$,
\begin{equation}
\|\mathbf R(\mathbf Y,C,s)\|
\le K e^{-s}\operatorname{diam}(\mathbf Y)
\le 2K e^{-s}\|\mathbf Y\|.
\label{eq:boundary-remainder-s-estimate}
\end{equation}
Moreover,
\begin{equation}
\begin{split}
&\|\widetilde{\mathbf R}(\mathbf Y_1,C_1,\eta)
-\widetilde{\mathbf R}(\mathbf Y_2,C_2,\eta)\|\\
&\qquad\le K\eta^2
\bigl(\|\mathbf Y_1-\mathbf Y_2\|+|C_1-C_2|\bigr)
\end{split}
\label{eq:boundary-remainder-Lipschitz}
\end{equation}
for all $(\mathbf Y_a,C_a)\in\mathcal U_Y\times\mathcal U_C$,
$a=1,2$.

The rescaled center term $F_C:=\eta^2\overline{U_c}$ also
admits a $C^1$ extension by zero to $\eta\le0$, with
\begin{equation}
|F_C|\le K\eta^2,
\qquad
\|D_{(\mathbf Y,C)}F_C\|\le K\eta^2,
\qquad
|\partial_\eta F_C|\le K\eta.
\label{eq:boundary-center-estimates}
\end{equation}
\end{proposition}

\begin{proof}
Write $b_*:=\operatorname{Im}C_*>0$. Choose the neighborhoods
so that $\operatorname{Im}C\ge b_*/2$ and
$\|\mathbf Y\|\le M$ throughout their product. Since
\[
D_{jm}
=C-\bar C+\eta(qY_j-\bar q\bar Y_m),
\]
we have
\[
|D_{jm}|\ge b_*-2M\eta.
\]
Thus, after reducing $\eta_0$ if necessary,
\begin{equation}
|D_{jm}|\ge b_*/2>0
\label{eq:boundary-denominator-bound}
\end{equation}
uniformly for all indices.

Set $S_j:=\sum_m\Gamma_m/D_{jm}$. The identity
$D_{jm}-D_{\ell m}=\eta q(Y_j-Y_\ell)$ gives
\[
S_j-S_\ell
=-\eta q(Y_j-Y_\ell)
\sum_{m=1}^3\frac{\Gamma_m}{D_{jm}D_{\ell m}}.
\]
Since
\[
U_j-U_c
=-\frac{1}{2\pi i\Gamma_{\mathrm{tot}}}
\sum_{\ell=1}^3\Gamma_\ell(S_j-S_\ell),
\]
we obtain the exact factorization
\begin{equation}
\begin{split}
R_j(\mathbf Y,C,-2\log\eta)
&=\eta^2q^2 A_j(\mathbf Y,C,\eta),\\
A_j(\mathbf Y,C,\eta)
&:=\frac{1}{2\pi i\Gamma_{\mathrm{tot}}}
\sum_{\ell,m=1}^3
\frac{\Gamma_\ell\Gamma_m(Y_j-Y_\ell)}
{D_{jm}D_{\ell m}}.
\end{split}
\label{eq:boundary-exact-factorization}
\end{equation}
By \eqref{eq:boundary-denominator-bound},
\[
\|\mathbf A\|\le K\operatorname{diam}(\mathbf Y).
\]
Complex conjugation preserves the norm, proving
\eqref{eq:boundary-remainder-size}.
Since $\eta^2=e^{-s}$, this also proves
\eqref{eq:boundary-remainder-s-estimate}.

The real derivatives of $\mathbf A$ with respect to
$(\mathbf Y,C)$ are uniformly bounded. Furthermore,
\[
\partial_\eta D_{jm}
=(1+2i\Omega)qY_j
-(1-2i\Omega)\bar q\bar Y_m
\]
is uniformly bounded. Together with
\eqref{eq:boundary-denominator-bound}, this implies that
$\partial_\eta\mathbf A$ is uniformly bounded.

Conjugating \eqref{eq:boundary-exact-factorization} yields
\[
\widetilde{\mathbf R}
=\eta^2\bar q^{\,2}\overline{\mathbf A}.
\]
Therefore,
\[
\partial_\eta\widetilde{\mathbf R}
=(2-4i\Omega)\eta\bar q^{\,2}\overline{\mathbf A}
+\eta^2\bar q^{\,2}\partial_\eta\overline{\mathbf A},
\]
which proves \eqref{eq:boundary-remainder-derivatives}.
Hence the function and all its first derivatives converge
uniformly to zero as $\eta\downarrow0$. The extension by zero
to $\eta\le0$ is therefore $C^1$, with derivative zero on
$\eta=0$. Applying the mean value theorem on the convex product
neighborhood gives \eqref{eq:boundary-remainder-Lipschitz}.

For the center term, direct conjugation gives
\[
F_C
=\frac{\eta^2}{2\pi i\Gamma_{\mathrm{tot}}}
\sum_{\ell,m=1}^3
\frac{\Gamma_\ell\Gamma_m}{\bar D_{\ell m}}.
\]
The same denominator and derivative bounds prove
\eqref{eq:boundary-center-estimates}, and hence the claimed
$C^1$ extension of $F_C$.

For completeness, the leading term can also be identified.
Let
\[
\mu:=\frac{1}{\Gamma_{\mathrm{tot}}}
\sum_{\ell=1}^3\Gamma_\ell Y_\ell.
\]
Using
\[
\frac{1}{D_{jm}D_{\ell m}}
=\frac{1}{(C-\bar C)^2}
+O(\eta\|\mathbf Y\|)
\]
and
\[
\sum_{\ell,m}\Gamma_\ell\Gamma_m(Y_j-Y_\ell)
=\Gamma_{\mathrm{tot}}^2(Y_j-\mu),
\]
we obtain
\begin{equation}
\begin{split}
\widetilde R_j
&=-\frac{\Gamma_{\mathrm{tot}}}
{2\pi i(C-\bar C)^2}
\eta^2e^{-4i\Omega\log\eta}
(\bar Y_j-\bar\mu)+\mathcal E_j,\\
\|\boldsymbol{\mathcal E}\|
&\le K\eta^3\|\mathbf Y\|^2.
\end{split}
\label{eq:boundary-leading-term}
\end{equation}
Thus the logarithmic rotation is compatible with the $C^1$
extension. If $\Omega\ne0$ and the components of $\mathbf Y$
are not all equal, however,
$\widetilde{\mathbf R}/\eta^2$ has no limit as
$\eta\downarrow0$, so a $C^2$ extension does not exist there.
\end{proof}

\medskip


\subsection{Existence of Interior Collapse via a Strong Stable Manifold}
\label{subsec:interior_collapse}



We establish the existence of finite-time collapsing solutions
for the perturbed three-vortex system by exploiting the stable
direction introduced by the auxiliary scaling variable.
An important feature of the argument is that the unperturbed
collapsing profile need not possess any stable eigenvalues
in the shape variables.

The autonomous augmented system is given by
\begin{equation}
\label{eq:augmented_recalled}
\begin{aligned}
\frac{d\mathbf Y}{ds}
&=G(\mathbf Y)+\widetilde R(\mathbf Y,C,\eta),\\
\frac{dC}{ds}
&=\eta^2\mathcal V_{\mathrm{drift}}(C,\mathbf Y,\eta),\\
\frac{d\eta}{ds}
&=-\frac12\eta,
\end{aligned}
\end{equation}
where
\[
G(\mathbf Y)
=
\left(\frac12+i\Omega\right)\mathbf Y
+V_{\mathrm{sing}}(\mathbf Y),
\qquad
\widetilde R(\mathbf Y,C,\eta)
=O(\eta^2|\mathbf Y|).
\]
Here $\mathbf Y\in\mathbb C^3$, $C\in\mathbb C^+$,
and $\eta\geq0$.

Let $\mathbf W^*$ be a  collision-free, non-equilateral
solution of the profile equation
$G(\mathbf W^*)=0.$
For any prescribed $C_*\in\mathbb H^+$, the point
$p^*=(\mathbf W^*,C_*,0)$
is an equilibrium of \eqref{eq:augmented_recalled}.

As shown in Appendix~\ref{sec:spectrum}, the spectral analysis of the unperturbed three-vortex
similarity system gives, for its real linearization
$L=DG(\mathbf W^*)$,
\begin{equation}
\label{eq:collapse_shape_spectrum}
\sigma(L)
=
\left\{
0,0,1,1,
\frac12+i\Omega,\frac12-i\Omega
\right\}.
\end{equation}
In particular, $L$ has no eigenvalues with strictly
negative real part. Nevertheless, the full augmented
system possesses a stable direction associated with
the equation $\eta'=-\eta/2$. This permits the
construction of collapsing trajectories by a strong
stable manifold argument.

\begin{theorem}[(Existence of Interior Collapse via a Strong Stable Manifold)]
\label{thm:interior_collapse}

Let $\Gamma_k\in\mathbb R\setminus\{0\}$,
$k=1,2,3$, satisfy
\[
\sum_{1\leq i<j\leq3}\Gamma_i\Gamma_j=0,
\]
and let $\mathbf W^*$ be a collision-free,
non-equilateral solution of the profile
equation $G(\mathbf W^*)=0$.

Then, for every prescribed collapse point
$C_*\in\mathbb H^+$, there exists a
one-dimensional local strong stable manifold
$\mathcal W^{ss}_{\mathrm{loc}}(p^*)$
of the augmented system, whose positive-$\eta$
branch is nonempty.

More precisely, there exist $\delta>0$ and
$C^1$ functions
\[
h_Y:[0,\delta)\to\mathbb C^3,
\qquad
h_C:[0,\delta)\to\mathbb C,
\]
satisfying
\[
h_Y(0)=h_C(0)=0,
\qquad
h_Y'(0)=h_C'(0)=0,
\]
such that the positive-$\eta$ branch is given locally by
\begin{equation}
\label{eq:stable_manifold_graph}
\mathcal W^{ss}_{\mathrm{loc},+}(p^*)
=
\left\{
\bigl(
\mathbf W^*+h_Y(\eta),
C_*+h_C(\eta),
\eta
\bigr):
0\leq\eta<\delta
\right\}.
\end{equation}
Every trajectory on this branch with $\eta>0$
satisfies
\begin{equation}
\label{eq:stable_asymptotics}
\mathbf Y(s)=\mathbf W^*+O(\eta(s)^2),
\qquad
C(s)=C_*+O(\eta(s)^2),
\qquad s\to+\infty.
\end{equation}

Suppose, in addition, that the augmented system
is related to the original vortex dynamics by the
exact similarity transformation
\begin{equation}
\label{eq:physical_reconstruction_collapse}
z_k(t)
=
C(s)+e^{-(1/2+i\Omega)s}Y_k(s),
\qquad
T-t=e^{-s}=\eta(s)^2.
\end{equation}
Then every trajectory on the positive-$\eta$ branch
generates an exact finite-time interior collapsing
solution satisfying
\begin{equation}
\label{eq:physical_collapse_rate}
z_k(t)
=
C_*+
\sqrt{T-t}\,e^{-i\Omega s(t)}W_k^*
+O(T-t),
\end{equation}
and, in particular,
\begin{equation}
\label{eq:collapse_rate_final}
|z_k(t)-C_*|
=
O\!\left(\sqrt{T-t}\right),
\qquad t\to T^-.
\end{equation}
The three vortices remain distinct and inside
$\mathbb C^+$ for all $t<T$ sufficiently close
to $T$, and collide simultaneously at $C_*$.
\end{theorem}

\begin{proof}
We divide the proof into three steps.

\medskip
\noindent
\textbf{Step 1. Linearization and the stable direction.}

Since $G(\mathbf W^*)=0$ and the perturbation
satisfies \eqref{eq:boundary-remainder-size},
the real Jacobian of the augmented vector field
at $p^*$ is
\begin{equation}
\label{eq:full_augmented_jacobian}
D\mathcal F(p^*)
=
\begin{pmatrix}
L&0&0\\
0&0_{2\times2}&0\\
0&0&-\frac12
\end{pmatrix}.
\end{equation}
Here all derivatives are understood in real
coordinates.

Indeed, the $C^1$ regularity and the
$O(\eta^2)$ estimate imply that the first
derivatives of $\widetilde R$ vanish at
$\eta=0$. The same holds for the center equation
because of its prefactor $\eta^2$.

It follows from \eqref{eq:collapse_shape_spectrum}
that
\[
\sigma(D\mathcal F(p^*))
=
\left\{
-\frac12,\,
0,0,0,0,\,
1,1,\,
\frac12+i\Omega,\frac12-i\Omega
\right\}.
\]
Consequently, the real phase space admits
an invariant spectral decomposition
\[
\mathbb R^9=E^s\oplus E^c\oplus E^u,
\]
where
\[
\dim E^s=1,\qquad
\dim E^c=4,\qquad
\dim E^u=4.
\]
The stable subspace is precisely
\[
E^s
=
\operatorname{span}_{\mathbb R}
\{(0,0,1)\}.
\]

Although the equilibrium is nonhyperbolic,
its stable eigenvalue $-1/2$ is separated
from the remaining spectrum. The strong
stable manifold theorem therefore gives a
one-dimensional local invariant manifold
tangent to $E^s$ at $p^*$.

In particular, the manifold can locally be
represented as a graph over the $\eta$-axis.
We now construct this graph and establish
the sharper estimates \eqref{eq:stable_asymptotics}.

\medskip
\noindent
\textbf{Step 2. Construction by a Lyapunov--Perron equation.}

Introduce the perturbation variables
\[
\mathbf y=\mathbf Y-\mathbf W^*,
\qquad
c=C-C_*.
\]
Since $\mathbf W^*$ is collision-free,
$G$ is smooth in a neighborhood of
$\mathbf W^*$. Hence
\[
G(\mathbf W^*+\mathbf y)
=
L\mathbf y+N(\mathbf y),
\qquad
N(\mathbf y)=O(|\mathbf y|^2).
\]
The augmented equations become
\begin{equation}
\label{eq:perturbed_augmented_equations}
\begin{aligned}
\mathbf y'
&=L\mathbf y+N(\mathbf y)
+\widetilde R(\mathbf W^*+\mathbf y,
             C_*+c,\eta),\\
c'
&=\eta^2\mathcal V_{\mathrm{drift}}
(C_*+c,\mathbf W^*+\mathbf y,\eta),\\
\eta'
&=-\frac12\eta.
\end{aligned}
\end{equation}

Fix a sufficiently small $\eta_0>0$ and
translate the initial similarity time to
$\tau=0$. The auxiliary equation gives
\begin{equation}
\label{eq:eta_solution_ss}
\eta(\tau)=\eta_0e^{-\tau/2},
\qquad \tau\geq0.
\end{equation}

We seek solutions for which
$\mathbf y(\tau)\to0$ and $c(\tau)\to0$
as $\tau\to+\infty$.

Because every eigenvalue of $L$ has
nonnegative real part, there exist constants
$M>0$ and $m\geq0$ such that
\begin{equation}
\label{eq:backward_semigroup_bound}
\|e^{-Lr}\|
\leq M(1+r)^m,
\qquad r\geq0.
\end{equation}
This estimate allows for polynomial growth
caused by possible Jordan blocks associated
with the zero eigenvalues.

Imposing the terminal conditions
$\mathbf y(\infty)=0$ and $c(\infty)=0$,
the variation-of-constants formula gives
the coupled integral equations
\begin{equation}
\label{eq:LP_shape}
\begin{aligned}
\mathbf y(\tau)
=-\int_\tau^\infty
&e^{L(\tau-u)}
\Bigl[
N(\mathbf y(u))\\
&\quad+
\widetilde R\bigl(
\mathbf W^*+\mathbf y(u),
C_*+c(u),
\eta_0e^{-u/2}
\bigr)
\Bigr]\,du
\end{aligned}
\end{equation}
and
\begin{equation}
\label{eq:LP_center}
\begin{aligned}
c(\tau)
=-\int_\tau^\infty
\eta_0^2e^{-u}\,
\mathcal V_{\mathrm{drift}}\bigl(
&C_*+c(u),\\
&\mathbf W^*+\mathbf y(u),
\eta_0e^{-u/2}
\bigr)\,du.
\end{aligned}
\end{equation}

Consider the Banach space
\[
\mathcal X
=
\left\{
(\mathbf y,c)\in
C([0,\infty);\mathbb R^6\times\mathbb R^2):
\|(\mathbf y,c)\|_{\mathcal X}<\infty
\right\},
\]
equipped with the weighted norm
\[
\|(\mathbf y,c)\|_{\mathcal X}
=
\sup_{\tau\geq0}
e^\tau\bigl(
|\mathbf y(\tau)|+|c(\tau)|
\bigr).
\]

By \eqref{eq:eta_solution_ss} and
\eqref{eq:boundary-remainder-s-estimate},
the forcing satisfies
\[
\left|
\widetilde R(
\mathbf W^*+\mathbf y,
C_*+c,\eta(\tau))
\right|
=O(\eta_0^2e^{-\tau})
\]
uniformly for $(\mathbf y,c)$ in a
sufficiently small neighborhood of the origin.
Similarly, the center drift is
$O(\eta_0^2e^{-\tau})$.

The semigroup estimate
\eqref{eq:backward_semigroup_bound} implies
\[
\int_\tau^\infty
\|e^{L(\tau-u)}\|e^{-u}\,du
\leq
M e^{-\tau}
\int_0^\infty(1+r)^m e^{-r}\,dr.
\]
The integral on the right is finite.

Together with $N(\mathbf y)=O(|\mathbf y|^2)$,
these bounds show that the right-hand sides
of \eqref{eq:LP_shape}--\eqref{eq:LP_center}
map a ball
\[
\|(\mathbf y,c)\|_{\mathcal X}
\leq K_0\eta_0^2
\]
into itself, for a suitable fixed $K_0>0$
and sufficiently small $\eta_0$.

Moreover, the local $C^1$ regularity proved in Proposition~\ref{prop:boundary-remainder-C1} implies
\[
D_{(\mathbf Y,C)}
\widetilde R(\mathbf Y,C,0)=0.
\]
Hence the Lipschitz constant of the remainder
with respect to $(\mathbf Y,C)$ tends to zero
as $\eta_0\to0$. The Lipschitz contributions
from $N$ and the center drift are likewise
small on the above ball.

It follows that the coupled integral operator
is a contraction for sufficiently small
$\eta_0$. By the Banach fixed-point theorem,
there exists a unique solution
$(\mathbf y,c)\in\mathcal X$ satisfying
\eqref{eq:LP_shape}--\eqref{eq:LP_center},
and
\begin{equation}
\label{eq:decay_estimate_ss}
|\mathbf y(\tau)|+|c(\tau)|
\leq K_0\eta_0^2e^{-\tau}.
\end{equation}

Define
\[
h_Y(\eta_0)=\mathbf y(0),
\qquad
h_C(\eta_0)=c(0).
\]
The parameter-dependent contraction argument,
together with the $C^1$ regularity of the
vector field, gives $C^1$ dependence on
$\eta_0$. Furthermore,
\eqref{eq:decay_estimate_ss} implies
\[
h_Y(\eta)=O(\eta^2),
\qquad
h_C(\eta)=O(\eta^2),
\]
and therefore
\[
h_Y'(0)=h_C'(0)=0.
\]

Uniqueness of the fixed point and autonomy
of the augmented system show that the graph
\eqref{eq:stable_manifold_graph} is forward
invariant. It is tangent to $E^s$ at $p^*$,
and thus coincides locally with the positive
branch of the strong stable manifold.

Since $\eta_0>0$ may be chosen arbitrarily
small, this branch is nonempty. Returning
to the original similarity time yields
\[
\mathbf Y(s)-\mathbf W^*=O(\eta(s)^2),
\qquad
C(s)-C_*=O(\eta(s)^2),
\]
which proves \eqref{eq:stable_asymptotics}.

\medskip
\noindent
\textbf{Step 3. Reconstruction and interior collapse.}

By the exact similarity transformation
\eqref{eq:physical_reconstruction_collapse},
\[
z_k(t)-C_*
=
C(s)-C_*
+
\eta(s)e^{-i\Omega s}Y_k(s).
\]
Using \eqref{eq:stable_asymptotics}, we find
\[
z_k(t)-C_*
=
\eta(s)e^{-i\Omega s}W_k^*
+O(\eta(s)^2).
\]
Since $\eta(s)^2=T-t$, this gives
\[
z_k(t)
=
C_*+
\sqrt{T-t}\,e^{-i\Omega s(t)}W_k^*
+O(T-t),
\]
and consequently
\[
|z_k(t)-C_*|
=
O\!\left(\sqrt{T-t}\right).
\]

Moreover, for $i\neq j$,
\[
z_i(t)-z_j(t)
=
\eta(s)e^{-i\Omega s}
\left(
W_i^*-W_j^*+O(\eta(s)^2)
\right).
\]
As the profile is collision-free,
$W_i^*\neq W_j^*$. Therefore,
\[
|z_i(t)-z_j(t)|
=
\sqrt{T-t}
\left(
|W_i^*-W_j^*|+O(T-t)
\right).
\]
Thus the vortices remain distinct for
$t<T$ sufficiently close to $T$, and
all three positions converge to $C_*$ as
$t\to T^-$.

Finally, because $C_*\in\mathbb C^+$,
we have $\operatorname{Im}C_*>0$.
By taking $\eta_0$ sufficiently small,
the preceding estimates ensure
\[
|z_k(t)-C_*|
<
\frac12\operatorname{Im}C_*
\]
throughout the constructed trajectory.
Hence
\[
\operatorname{Im}z_k(t)>0,
\qquad k=1,2,3,
\]
and the collapse takes place strictly in
the interior of the upper half-plane.

This completes the proof.
\end{proof}

\begin{remark}
These collapse solutions are asymptotically self-similar, which means that the collapsing vortices approach a fixed configuration when you zoom in and rotate with them.
Let \(z_j(t)\in\mathbb H^+\) be the vortex positions, \(C_*\) the collapse point, and \(T\) the collapse time. The usual square-root, logarithmically rotating form is
\[
\frac{e^{i\Omega\log(T-t)}}{\sqrt{T-t}}\bigl(z_j(t)-C_*\bigr)
\;\longrightarrow\; W_j^*
\qquad\text{as }t\uparrow T,
\]
where \(\Omega\) is a real constant and the distinct points \(W_j^*\) form the limiting profile. In particular, each pairwise distance satisfies
\[
|z_i(t)-z_j(t)|\sim |W_i^*-W_j^*|\sqrt{T-t}.
\]
For three vortices, this means the triangle's side-length ratios converge even as the triangle shrinks and can keep rotating.
\end{remark}

\begin{remark}[(Existence versus stability)]
\label{rem:collapse_not_stable}
The theorem establishes the existence of
collapsing trajectories, not Lyapunov stability
of the collapsing profile. The augmented
equilibrium has positive-real-part eigenvalues,
and its local strong stable manifold is only
one-dimensional in the unreduced nine-dimensional
phase space. Thus the construction selects
special initial configurations for which the
growing shape modes are suppressed.
Generic nearby initial conditions need not
converge to the same collapsing profile.

\end{remark}


\medskip


\begin{appendices}

\section{The complete spectrum of the real \(6\times6\) Jacobian}\label{sec:spectrum}

Under the assumptions that \(\Gamma_1,\Gamma_2,\Gamma_3\in\mathbb R\setminus\{0\}\), \(\Omega\in\mathbb R\), and \(\mathbf Y^*\) is a collision-free equilibrium, we will show that
\[
\sigma\big(DG(\mathbf Y^*)\big)
=
\left\{
0,0,1,1,
\frac12+i\Omega,
\frac12-i\Omega
\right\}.
\]
This is a consequence of the special Hamiltonian structure and symmetries of the three-vortex problem, rather than an approximation. 

\subsection{Derivation of the real \(6\times6\) Jacobian}
Consider the system
\begin{equation}\label{eq:Y_G}
\frac{dY_k}{ds}
=G_k(\mathbf Y)
=
\left(\frac12+i\Omega\right)Y_k
-\frac{1}{2\pi i}
\sum_{j\ne k}^{3}
\frac{\Gamma_j}{\bar Y_k-\bar Y_j}.
\end{equation}
Let \(\mathbf Y^*\) be an equilibrium, i.e., $G(\mathbf Y^*)=0.$
We also introduce a perturbation
\[
\mathbf Y(s)=\mathbf Y^*+\boldsymbol\eta(s),
\qquad
\boldsymbol\eta\in\mathbb C^3,
\]
and define $a=\frac12+i\Omega.$
Expanding the interaction term to first order gives
\[
\frac{1}{\bar Y_k-\bar Y_j}
=
\frac{1}{\bar Y_k^*-\bar Y_j^*}
-
\frac{\bar\eta_k-\bar\eta_j}
{(\bar Y_k^*-\bar Y_j^*)^2}
+O(\|\boldsymbol\eta\|^2).
\]
Therefore, the linearized system is
\begin{equation}\label{eq:linear_YG}
\frac{d\boldsymbol\eta}{ds}
=
a\boldsymbol\eta+B\bar{\boldsymbol\eta},
\end{equation}
where \(B\in\mathbb C^{3\times3}\) has entries
\begin{equation}\label{eq:B}
B_{kj}=
\begin{cases}
\displaystyle
\frac{1}{2\pi i}
\sum_{\ell\neq k}
\frac{\Gamma_\ell}
{(\bar Y_k^*-\bar Y_\ell^*)^2},
&j=k,\\[12pt]
\displaystyle
-\frac{\Gamma_j}
{2\pi i(\bar Y_k^*-\bar Y_j^*)^2},
&j\neq k.
\end{cases}
\end{equation}
Notice that the perturbation involves both \(\boldsymbol\eta\) and its complex conjugate. Consequently, \(G\) is not holomorphic, and its linearization is real-linear rather than complex-linear.

\subsubsection*{Real representation}

Write $\boldsymbol\eta=\mathbf u+i\mathbf v, \,B=P+iQ,$
where \(\mathbf u,\mathbf v\in\mathbb R^3\) and \(P,Q\in\mathbb R^{3\times3}\).
Then, one can rewrite \eqref{eq:linear_YG} as follows.
\[
\begin{aligned}
\mathbf u'
&=\left(\frac12 I_3+P\right)\mathbf u
+\left(-\Omega I_3+Q\right)\mathbf v,\\
\mathbf v'
&=\left(\Omega I_3+Q\right)\mathbf u
+\left(\frac12 I_3-P\right)\mathbf v.
\end{aligned}
\]
Thus, for
$
\boldsymbol\xi=
\begin{pmatrix}
\mathbf u\\
\mathbf v
\end{pmatrix}\in\mathbb R^6,
$
we obtain
\begin{equation}\label{eq:linear_real}
\frac{d\boldsymbol\xi}{ds}
=J\boldsymbol\xi,
\qquad
J=
\begin{pmatrix}
\frac12 I_3+P&-\Omega I_3+Q\\
\Omega I_3+Q&\frac12 I_3-P
\end{pmatrix}.
\end{equation}
This is precisely the real Jacobian \(J=DG(\mathbf Y^*)\) that we want to study.

An immediate observation is
\[
\operatorname{tr}J=3,
\]
independently of the equilibrium and the circulations.

\subsection{The eigenvalues obtained directly from symmetries}

\subsubsection*{Translation modes: \(\lambda=\frac12\pm i\Omega\)}

The vortex interaction depends only on differences \(Y_k-Y_j\).
Consequently,
\[
B\mathbf e=0,
\qquad
\mathbf e=(1,1,1)^T.
\]
Consider perturbations corresponding to translating all three vortices by the same complex number:
$\boldsymbol\eta(s)=c(s)\mathbf e.$
Equation \eqref{eq:linear_YG} reduces to
\[
c'=\left(\frac12+i\Omega\right)c.
\]
The two eigenvalues of this real two-dimensional translation subspace are therefore
\begin{equation}\label{eq:trans_eigen}
\lambda_{T,\pm}=\frac12\pm i\Omega.
\end{equation}

\subsubsection*{Rotation mode: \(\lambda=0\)}

The vortex vector field is rotationally equivariant:
$$
G(e^{i\theta}\mathbf Y)
=
e^{i\theta}G(\mathbf Y).
$$
Since \(G(\mathbf Y^*)=0\), every rotated configuration \(e^{i\theta}\mathbf Y^*\) is also an equilibrium.
Differentiating with respect to \(\theta\) at zero gives
\[
DG(\mathbf Y^*)[i\mathbf Y^*]=0.
\]
Hence,
\begin{equation}\label{eq:rot_eigen}
\lambda_R=0,
\qquad
\boldsymbol\eta_R=i\mathbf Y^*.
\end{equation}
This is a neutral eigenmode associated with rotational symmetry.

\subsubsection*{Scaling mode: \(\lambda=1\)}

Define the original vortex interaction
\[
V(\mathbf Y)
=
-\frac{1}{2\pi i}
\sum_{j\ne k}^3
\frac{\Gamma_j}{\bar Y_k-\bar Y_j}.
\]
Then,
$
G(\mathbf Y)=a\mathbf Y+V(\mathbf Y).
$
Because the vortex interaction is homogeneous of degree \(-1\),
\[
V(r\mathbf Y)=r^{-1}V(\mathbf Y),
\qquad r>0.
\]
Differentiating at \(r=1\),
\[
DV(\mathbf Y^*)[\mathbf Y^*]=-V(\mathbf Y^*).
\]
The equilibrium condition gives
$
V(\mathbf Y^*)=-a\mathbf Y^*.
$
Consequently,
\begin{equation}\label{eq:scale_eigen_comp}
DG(\mathbf Y^*)[\mathbf Y^*]
=(1+2i\Omega)\mathbf Y^*.
\end{equation}
This does not mean that \(1+2i\Omega\) is an eigenvalue of the real Jacobian, because the right-hand side is not generally a real scalar multiple of \(\mathbf Y^*\).

Instead, combining \eqref{eq:rot_eigen} and \eqref{eq:scale_eigen_comp}, we have
\[
DG(\mathbf Y^*)
[\mathbf Y^*+2\Omega i\mathbf Y^*]
=
\mathbf Y^*+2\Omega i\mathbf Y^*.
\]
Thus,
\begin{equation}\label{eq:scale_eigen}
\lambda_S=1,
\qquad
\boldsymbol\eta_S=(1+2i\Omega)\mathbf Y^*.
\end{equation}
This growing mode is associated with scaling and, physically, a perturbation of the collapse time.

At this stage, we have identified four of the six eigenvalues:
$
\frac12+i\Omega,\;
\frac12-i\Omega,\;
0,\; 1.
$
The remaining two require exploiting the special Hamiltonian structure of the three-vortex dynamics.

\subsection{Determining the remaining two eigenvalues}

The key is to use the Hamiltonian structure of the original point-vortex dynamics.

\subsubsection{Necessary conditions for a collapsing equilibrium}

Introduce the circulation-weighted quantities
\[
C=\sum_{k=1}^3\Gamma_kY_k^*,
\qquad
I_\Gamma=\sum_{k=1}^3\Gamma_k|Y_k^*|^2,
\]
and
\[
L_\Gamma
=\Gamma_1\Gamma_2+\Gamma_2\Gamma_3+\Gamma_3\Gamma_1.
\]
The necessary conditions for nontrivial three-vortex self-similar collapse gives
$$
C=0,\qquad I_\Gamma=0,
\qquad
L_\Gamma=0.
$$

In particular,
$
\left(\sum_k\Gamma_k\right)^2
=\sum_k\Gamma_k^2>0,
$
so the total circulation is nonzero.

\subsubsection{The centered perturbation space}

Define
\[
\mathcal V
=
\left\{
\boldsymbol\eta\in\mathbb C^3:
\sum_{k=1}^3\Gamma_k\eta_k=0
\right\}.
\]
As a real vector space, \(\mathcal V\) has dimension four.
The real Jacobian preserves \(\mathcal V\), because
\[
\sum_{k=1}^3\Gamma_kG_k(\mathbf Y)
=
a\sum_{k=1}^3\Gamma_kY_k.
\]
Moreover, the full perturbation space splits into two invariant subspaces,
$
\mathbb R^6=\mathcal T\oplus\mathcal V,
$
here \(\mathcal T\) is the two-dimensional translation subspace.
Thus, the remaining four eigenvalues belong to \(J|_{\mathcal V}\).

\subsubsection{The conformally symplectic identity}

Let
$
D_\Gamma=\operatorname{diag}(\Gamma_1,\Gamma_2,\Gamma_3),
$
and introduce the skew-symmetric matrix
\[
W=
\begin{pmatrix}
0&D_\Gamma\\
-D_\Gamma&0
\end{pmatrix}.
\]
The original point-vortex interaction is Hamiltonian with respect to this circulation-weighted symplectic structure.

The rotation term \(i\Omega\mathbf Y\) is also Hamiltonian, whereas the scaling term \(\frac12\mathbf Y\) gives a conformal contribution.
Consequently, the real Jacobian $J=DG$ satisfies
\begin{equation}\label{eq:J_W}
J^TW+WJ=W.
\end{equation}

\begin{proof}[Proof of Equation \eqref{eq:J_W}]

We can write
\[
G(\mathbf Y)
=
\underbrace{V(\mathbf Y)}_{\text{vortex interaction}}
+
\underbrace{i\Omega\mathbf Y}_{\text{rotation}}
+
\underbrace{\frac12\mathbf Y}_{\text{scaling}}.
\]
Next we investigate the three terms separately.

Write \(Y_k=x_k+iy_k\), and introduce the circulation-weighted symplectic form
$
\omega_\Gamma
=
\sum_{k=1}^3\Gamma_k\,dx_k\wedge dy_k.
$
The corresponding matrix is  $W$ defined above
\[
W=
\begin{pmatrix}
0&D_\Gamma\\
-D_\Gamma&0
\end{pmatrix},
\qquad
D_\Gamma=\operatorname{diag}(\Gamma_1,\Gamma_2,\Gamma_3).
\]
A vector field \(X\) is Hamiltonian if it can be generated by a Hamiltonian function \(H\) through
$
\iota_X\omega_\Gamma=dH.
$
For a Hamiltonian vector field, its Jacobian \(A=DX\) satisfies
\[
A^TW+WA=0.
\]
Geometrically, its flow preserves \(\omega_\Gamma\).

The original vortex interaction \(V\) is Hamiltonian, with
$
H_V
=
-\frac{1}{2\pi}
\sum_{j<k}\Gamma_j\Gamma_k
\log|Y_j-Y_k|.
$
Therefore,
\begin{equation}\label{eq:DV_W}
DV(\mathbf Y)^TW+WDV(\mathbf Y)=0.
\end{equation}

Consider $\mathbf Y'=i\Omega\mathbf Y.$
For each vortex,
$Y_k'=i\Omega Y_k.$
Writing \(Y_k=x_k+iy_k\),
\[
\begin{pmatrix}
x_k'\\y_k'
\end{pmatrix}
=
\begin{pmatrix}
0&-\Omega\\
\Omega&0
\end{pmatrix}
\begin{pmatrix}
x_k\\y_k
\end{pmatrix}.
\]
Thus, in \(\mathbb R^6\), the Jacobian of the rotation field is
\[
R=
\begin{pmatrix}
0&-\Omega I_3\\
\Omega I_3&0
\end{pmatrix}.
\]
Direct matrix multiplication yields
\begin{equation}\label{eq:R_W}
R^TW+WR=0.
\end{equation}
Therefore, rotation is Hamiltonian.

The Jacobian of the scaling field is
$
S=\frac12I_6.
$
Therefore,
\[
S^TW+WS
=\frac12W+\frac12W
=W.
\]
Hence,
\begin{equation}\label{eq:S_W}
S^TW+WS=W.
\end{equation}
Compare this with the Hamiltonian condition \(A^TW+WA=0\). The nonzero right-hand side shows that scaling is not Hamiltonian with respect to the fixed symplectic form \(\omega_\Gamma\).

Finally, adding the identities \eqref{eq:DV_W}, \eqref{eq:R_W} and \eqref{eq:S_W} gives \eqref{eq:J_W}
\end{proof}

Identity \eqref{eq:J_W} implies that the spectrum is symmetric under
$
\lambda\longmapsto1-\lambda.
$
More precisely, if \(\lambda\) is an eigenvalue, then \(1-\lambda\) is also an eigenvalue, with the same algebraic multiplicity.

To see this,  from
$
J^TW+WJ=W,
$
and the invertibility of \(W\), we obtain
\[
W^{-1}J^TW=I_6-J.
\]
Thus, \(J^T\) is similar to \(I_6-J\).

Since \(J\) and \(J^T\) have the same characteristic polynomial,
\[
\lambda\in\sigma(J)
\quad\Longrightarrow\quad
1-\lambda\in\sigma(J).
\]
This gives a spectral symmetry about the vertical line
$
\operatorname{Re}\lambda=\frac12.
$

\subsubsection{Exploiting the zero angular impulse}

Define the weighted symplectic form
$
\omega_\Gamma(\mathbf u,\mathbf v)
=
\sum_{k=1}^3
\Gamma_k\operatorname{Im}(\bar u_kv_k).
$
Consider the real two-dimensional subspace
\[
\mathcal E=
\operatorname{span}_{\mathbb R}
\{\mathbf Y^*,i\mathbf Y^*\}.
\]
We have already shown that this subspace is invariant under \(J\), with eigenvalues \(0\) and \(1\).

Furthermore,
\[
\omega_\Gamma(\mathbf Y^*,i\mathbf Y^*)
=
\sum_k\Gamma_k|Y_k^*|^2
=0.
\]
Thus \(\mathcal E\) is isotropic, i.e., $\omega_\Gamma(\mathbf u,\mathbf v)=0\,\text{for all }\mathbf u,\mathbf v\in\mathcal E$. Since it has dimension two inside the four-dimensional symplectic space \(\mathcal V\), it is a Lagrangian subspace.

In suitable coordinates adapted to \(\mathcal E\), the centered Jacobian has block form
\begin{equation}\label{eq:J_V_block}
J|_{\mathcal V}
\sim
\begin{pmatrix}
A&E\\
0&I_2-A^T
\end{pmatrix},
\end{equation}
where \(A\) is the restriction of \(J\) to \(\mathcal E\).

\begin{proof}[Proof of \eqref{eq:J_V_block}]
The advantage of a Lagrangian subspace is that we can construct a particularly convenient basis.
Suppose
$
\mathcal E=\operatorname{span}\{e_1,e_2\}.
$
Because \(\mathcal E\) is Lagrangian, we can choose two additional vectors \(f_1,f_2\) so that
\[
\mathcal V
=
\operatorname{span}\{e_1,e_2,f_1,f_2\},
\]
and
\[
\begin{aligned}
\omega_\Gamma(e_i,e_j)&=0,\\
\omega_\Gamma(f_i,f_j)&=0,\\
\omega_\Gamma(e_i,f_j)&=\delta_{ij}.
\end{aligned}
\]
This is a symplectic basis, obtainable by symplectic Gram--Schmidt.

In these coordinates, every vector in \(\mathcal V\) can be written as
\[
\boldsymbol\eta
=
q_1e_1+q_2e_2+p_1f_1+p_2f_2.
\]
Define
$
q=
\begin{pmatrix}q_1\\q_2\end{pmatrix},
\;
p=
\begin{pmatrix}p_1\\p_2\end{pmatrix}.
$
Then
$
\boldsymbol\xi=
\begin{pmatrix}q\\p\end{pmatrix}\in\mathbb R^4.
$
In these coordinates the symplectic matrix becomes
\[
W_0=
\begin{pmatrix}
0&I_2\\
-I_2&0
\end{pmatrix}.
\]
The important advantage is that the subspace \(\mathcal E\) corresponds simply to
\[
\mathcal E=\{(q,p):p=0\}.
\]
Thus, the first two coordinates describe perturbations within the scale-and-rotation plane, whereas the remaining coordinates describe directions complementary to that plane.

Let
$
J_{\mathcal V}
=
DG(\mathbf Y^*)|_{\mathcal V}.
$
In our new coordinates, we know that \(\mathcal E\) is invariant under the Jacobian:
$
J_{\mathcal V}\mathcal E\subseteq\mathcal E.
$
Eevery vector in \(\mathcal E\) has \(p=0\).
Therefore,
\[
J_{\mathcal V}
=
\begin{pmatrix}
A&K\\
0&D
\end{pmatrix}.
\]
The zero lower-left block comes from invariance of \(\mathcal E\), not from the symplectic condition itself.

Now we use the conformally symplectic identity derived earlier:
$
J_{\mathcal V}^TW_0+W_0J_{\mathcal V}=W_0.
$
We obtain
\[
J_{\mathcal V}
=
\begin{pmatrix}
A&K\\
0&D
\end{pmatrix},
\qquad
W_0=
\begin{pmatrix}
0&I_2\\
-I_2&0
\end{pmatrix}.
\]
By calculating the two products explicitly, we get
$
D=I_2-A^T,
\; K=K^T.
$
Therefore,
\[
J_{\mathcal V}
\sim
\begin{pmatrix}
A&K\\
0&I_2-A^T
\end{pmatrix},
\qquad K=K^T.
\]
This is exactly the block structure stated in \eqref{eq:J_V_block}. The symbol \(\sim\) indicates that the original real Jacobian restricted to \(\mathcal V\) is similar to this matrix through a change of basis.
\end{proof}

Since
$
\sigma(A)=\{0,1\},
$
we also have
$
\sigma(I_2-A^T)=\{1,0\}.
$
Therefore,
\begin{equation}\label{eq:J_V_eigen}
\sigma(J|_{\mathcal V})=\{0,0,1,1\}.
\end{equation}
This determines the remaining eigenvalues.

\subsection{Complete characteristic polynomial and spectrum}

Combining the translation and centered modes gives
\[
\det(\lambda I_6-J)
=
\lambda^2(\lambda-1)^2
\left[
\left(\lambda-\frac12\right)^2+\Omega^2
\right].
\]
Hence the six eigenvalues, counted with algebraic multiplicity, are
\[
\lambda_1=\lambda_2=0,\quad
\lambda_3=\lambda_4=1,\quad
\lambda_{5,6}=\frac12\pm i\Omega.
\]
The real parts of the eigenvalues are therefore
$
0,\, 0,\, 1,\,1,\, \frac12,\, \frac12.
$
Interestingly, these real parts do not depend on the particular circulations or the equilibrium shape, provided the assumptions hold.

Since some eigenvalues have strictly positive real parts, the equilibrium is linearly unstable in the full, unreduced similarity system.
An important subtlety is that the eigenvalues \(0\) and \(1\) each have algebraic multiplicity two. The Jacobian need not be diagonalizable. In particular, Jordan blocks of size two can occur, producing terms of the form
$
s e^{\lambda s}
$
in the linearized evolution.
Thus, knowing the eigenvalues alone does not completely determine perturbation growth.

\end{appendices}

\medskip


\end{document}